\documentclass[11pt, reqno]{amsart}
\usepackage{amscd}
\usepackage{latexsym}
\usepackage{amssymb}
\usepackage{amsmath}
\usepackage{amsthm}
\usepackage{indentfirst}
\usepackage{pifont}
\usepackage{graphicx}
\usepackage{geometry}
\usepackage{subfig}
\usepackage{cases}
\usepackage{hyperref}
\usepackage{cite}
\usepackage{color}

\numberwithin{equation}{section}

\newtheorem{theorem}{Theorem}[section]
\newtheorem{lemma}{Lemma}[section]

\theoremstyle{definition}

\newtheorem{proposition}{Proposition}[section]

\theoremstyle{remark}

\newcommand{\di}{\displaystyle}
\newcommand{\x}{\xi}
\newcommand{\RR}{\mathbb{R}}
\newcommand{\mb}{\mathbf}
\newcommand{\ra}{\rightarrow}

\newcommand{\ve}{\varepsilon}

\begin{document}


\title[Stability of degenerate Oleinik shock and rarefaction waves]{Time-asymptotic stability of composite waves of degenerate Oleinik shock and rarefaction for non-convex conservation laws}

\author[Huang]{Feimin Huang}
\address[Feimin Huang]{\newline Institute of Applied Mathematics, AMSS, CAS, Beijing 100190, P. R. China
\newline
and School of Mathematical Sciences, University of Chinese Academy of Sciences,
\newline Beijing 100049, P. R. China}
\email{fhuang@amt.ac.cn}

\author[Wang]{Yi Wang}
\address[Yi Wang]{\newline Institute of Applied Mathematics, AMSS, CAS, Beijing 100190, P. R. China
\newline
and School of Mathematical Sciences, University of Chinese Academy of Sciences,
\newline Beijing 100049, P. R. China}
\email{wangyi@amss.ac.cn}

\author[Zhang]{Jian Zhang}
\address[Jian Zhang]{\newline Institute of Applied Mathematics, AMSS, CAS, Beijing 100190, P. R. China
\newline
and School of Mathematical Sciences, University of Chinese Academy of Sciences,
\newline Beijing 100049, P. R. China}
\email{zhangjian2020@amss.ac.cn}

\date{}
\maketitle

\vspace{-0.6cm}

\noindent{\bf Abstract}: 
We are concerned with the large-time behavior of  the solution to one-dimensional (1D) cubic non-convex scalar viscous conservation laws. Due to the inflection point of the cubic non-convex flux,  the solution to the corresponding inviscid Riemann problem  can be the composite wave of a degenerate Oleinik shock and a rarefaction wave and these two nonlinear waves are always attached together. We give a first proof of the time-asymptotic stability of this composite wave, up to a time-dependent shift to the viscous Oleinik shock,  for the viscous equation. The Oleinik shock wave strength can be arbitrarily large. The main difficulty is due to the incompatibility of the time-asymptotic stability proof framework of individual viscous shock by the so-called anti-derivative method and the direct $L^2$-energy method to rarefaction wave. Here we develop a new type of $a$-contraction method with suitable weight function  and  the time-dependent shift to the viscous shock, which is motivated by \cite{KV,  KVW}. Another difficulty comes from that the Oleinik shock and rarefaction wave are always attached together and their wave interactions are very subtle. Therefore, the same time-dependent shift needs to be equipped to both Oleinik shock and rarefaction wave such that the wave interactions can be treated in our stability proof. Time-asymptotically, this shift function grows strictly sub-linear with respect to the time and then the shifted rarefaction wave is equivalent to the original self-similar rarefaction wave.



 \section{Introduction}
 We are concerned with the large-time behavior of the solution to the following one-dimensional (1D)  cubic non-convex scalar conservation laws
 \begin{equation}\label{1.1}
\left\{\begin{array}{ll}
u_t+f(u)_x=\mu u_{xx},&\di \quad f(u)=u^3,\quad (t,x)\in  \mathbb{R}_{+}\times \mathbb{R},\\[2mm]
u(0,x)=u_0(x)\ra u_\pm,&\di \quad {\rm as}\ x\ra \pm\infty, \\[2mm] 
\end{array}\right.
\end{equation}
where $u=u(t,x)\in \RR$ is the unknown function, the so-called conserved quantity and $\mu>0$ is the viscosity coefficient, $u_0(x)$ is the given initial data, and $u_{\pm}\in\mathbb{R}$ are the prescribed far field states.  
Equation \eqref {1.1} is closely related to rotationally invariant hyperbolic waves in viscoelasticity and magnetohydrodynamics \cite{BHR, LTP3} and can also be derived from nonlinear electromagnetic waves \cite{IPD,GAN,MEV}.

 The large-time asymptotic behavior of the solutions to  \eqref{1.1} with different far fields  $u_{\pm}$ is expected to be determined  by the corresponding inviscid Riemann problem 
\begin{equation}\label{1.2}
\left\{\begin{array}{l}
u_t+f(u)_x=0,
\\
u(0,x)=u_0(x)=
\left\{\begin{array}{ll}
u_-,&\di x<0,\\
u_+,&\di x>0.
\end{array}\right.
\end{array}\right.
\end{equation}

When the flux $f(u)$ is strictly convex (e.g., the classical Burgers equation), the solution to the Riemann problem \eqref{1.2} is either a classical shock (for $u_->u_+$) or a rarefaction wave  (for $u_-<u_+$), and the time-asymptotic stability of the corresponding viscous shock wave and rarefaction wave to its viscous conservation laws are extensively studied and well-established since the pioneer work of Il'in-Oleinik \cite{IO}. However, when the flux $f(u)$ is non-convex with inflection points, the wave structure of Riemann problem \eqref{1.2} can be much more complicated and consist of multiple wave patterns including Oleinik shock, rarefaction waves and even contact discontinuities.  Here Oleinik shock means that the shock wave satisfies the Oleinik entropy condition due to the non-convexity of the flux.

 There are rich literatures on the time-asymptotic stability of a single viscous Oleinik shock wave with the non-convex  flux $f(u)$ for the scalar viscous conservation laws. When  the non-convex  flux  $f(u)$ has only one inflection point, Kawashima-Matsumura\cite{MK} first proved the stability of  non-degenerate Oleinik shock wave by the anti-derivative and $L^2$-weighted energy methods, and then Mei \cite{ME} proved  the stability of its degenerate Oleinik shock wave.  For the general non-convex  flux  $f(u)$, Matsumura-Nishihara\cite{MN} proved the stability of both non-degenerate and  degenerate Oleinik shock wave by using a suitably choosing unbounded weighted function. See also Liu \cite{LHL} for the stability of non-degenerate shock wave by a different weight function and Huang-Xu \cite{HX} for the decay rate toward the viscous shock profile. Furthermore, Freist$\ddot{\mbox{u}}$hler-Serre \cite{FS} established an interesting $L^1$ stability theorem and Jones-Gardner-Kapitula \cite{JGK}, Osher-Ralston \cite{OR} and Weinberger \cite{We} also investigated the time-asymptotic stability of viscous shock wave by different methods. 
 
 For a single rarefaction wave, it can only be formed on the convex (or concave) part of the flux, therefore, the stability of the single rarefaction wave is exactly same as the strictly convex  (or concave)  case by using the direct $L^2$-energy method. Therefore, the time-asymptotic stability of a single Oleinik shock (including both non-degenerate and degenerate cases) or a single rarefaction wave to the non-convex scalar viscous conservation laws is rather satisfactory. 
 
However, as we described before, the wave structure  of Riemann problem \eqref{1.2}  for the non-convex flux $f(u)$ can be much more complicated and consist of multiple wave patterns including both Oleinik shock and rarefaction wave.  And there is no result on the time-asymptotic stability of these multiple wave patterns including both Oleinik shock and rarefaction wave to the viscous equation as far as we know. In the present paper, we are concerned with the cubic non-convex flux with one inflection point (i. e., $f(u)=u^3$) and prove the time-asymptotic stability of the composite wave including both degenerate Oleinik shock and rarefaction wave to scalar viscous equation \eqref{1.1}.

Compared with the single wave pattern case, the main difficulty here is due to the incompatibility of the time-asymptotic stability proof framework of individual viscous shock by the so-called anti-derivative method and the direct $L^2$-energy method to rarefaction wave. Here we develop a new type of $a$-contraction method with suitable weight function  and  the time-dependent shift to the viscous shock, which is motivated by the recent works \cite{KV, KVW}. In \cite{KV}, Kang-Vasseur first developed the $a$-contraction method for viscous Burgers equation  and  proved $L^2$-contraction of viscous shock wave up to a time-dependent shift function.  Then this idea is used in \cite{VKM} to obtain $L^2$-contraction of small viscous shock wave to the 1D barotropic Navier-Stokes equations by choosing the suitable weight function. $L^2$-contraction of shock waves were also extensively studied in \cite{MJ, KVW2}. Note that the contraction property in \cite{VKM,MJ}  holds for small shock essentially due to the choice of weight function. Here we develop a new $a$-contraction with suitable weight function such that the Oleinik shock wave strength can be arbitrarily large. The main idea is explained in Section 2 and the details can be found in Section 4.

Another motivation is the recent work of Kang-Vasseur-Wang \cite{KVW} for the time-asymptotic stability of the composite wave of viscous shock and rarefaction wave to 1D compressible barotropic Navier-Stokes equations, see also \cite{KVW1} for the  time-asymptotic stability of generic Riemann solutions to the compressible Navier-Stokes-Fourier equations. However, the main difference here from \cite{KVW}  is that the Oleinik shock and rarefaction wave in \eqref{1.2} are always attached together, while the viscous shock and rarefaction waves for the compressible Navier-Stokes equations in \cite{KVW}  are separate in  two different genuinely nonlinear characteristic fields of the corresponding Euler system. Therefore, the wave interactions between degenerate viscous Oleinik shock and rarefaction wave are  more subtle here. Since the time-dependent shift function for the viscous Oleinik shock is essential in our stability proof, we need to equip the same time-dependent shift to rarefaction wave such that two shifted waves are still attached together, otherwise, the shifted shock wave will cross through the rarefaction wave  and then the wave interactions will cause substantial difficulties. Time-asymptotically, we can prove that the time-dependent shift function grows strictly sub-linear with respect to the time $t$ and then the shifted rarefaction wave is still equivalent to the original self-similar rarefaction wave.

 The rest part of this paper is organized as follows. In Section 2, we first list some properties of the viscous Oleinik shock and rarefaction wave, and then present the new ideas for the $a$-contraction method, and  finally state our main result on the time-asymptotic stability of composite wave of degenerate Oleinik shock and rarefaction wave. Section 3 is devoted to the proof  of our main result. An approximate rarefaction wave is first constructed and then the weight function and time-dependent shift are well established. Then our stability problem is reformulated to the time-asymptotic stability of the shifted Oleinik shock and shifted approximate rarefaction wave. Finally, the main result 
 can be proved by the continuity arguments based on the local existence of the solution and the uniform-in-time a priori estimates. In Section 4, we establish the uniform-in-time a priori estimates by using the modified $a$-contraction method and the subtle wave interaction estimates between the Oleinik shock and rarefaction wave. Section 5 is an appendix for some detailed calculations.
 
   \section{Preliminaries and main result} 
   In this section, we start with the descriptions of the viscous Oleinik shock wave and the rarefaction wave. Then we explain the main idea of our modified $a$-contraction method and state our main result on the time-asymptotic  stability of composite wave of degenerate viscous Oleinik shock and rarefaction for the viscous equation \eqref{1.1}.
   
\noindent{\bf Oleinik shock wave.}  A shock wave with speed $\sigma$ for the non-convex conservation law \eqref{1.2} with the Riemann initial data $u_\pm$ is admissible if it satisfies the Rankine-Hugoniot condition 
 \begin{equation}\label{2.1}
\sigma=\frac{f(u_+)-f(u_-)}{u_+-u_-},
\end{equation}
     and the  Oleinik (or Liu) entropy condition
\begin{equation}\label{2.2}
\sigma \leq\frac{f(u)-f(u_-)}{u-u_-},\ \ \mbox{for} \ \mbox{all} \ u\ \mbox{between} \ u_-  \ \mbox{and} \ u_+.
\end{equation}
One can refer to the recent monagraph \cite{LT} for details.

In the present paper, we are concerned with non-classical shock wave, that is the shock wave will cross through the inflection point $u=0$ for the non-convex flux $f(u)=u^3$. From \eqref{2.1} and \eqref{2.2},  a shock wave for \eqref{1.2}  
with $f(u)=u^3$ is admissible, or called Oleinik shock, if for $u_-<0,$ the right state $u_+$ satisfies that
 \begin{equation}\label{2.3}
 u_-<u_+\leq-\frac{u_-}{2},
 \end{equation}
 or for $u_->0,$
 \begin{equation}\label{2.4}
- \frac{u_-}{2}\leq u_+ < u_-.
 \end{equation}
 The degenerate case $u_+=- \frac{u_-}{2}$ is calculated through the relation
 $$
\sigma=f^\prime(u_+),
 $$
 which means that the shock speed $\sigma$ is equal to the right characteristic speed $f^\prime(u_+)$.
In particular, when $u_+=-\frac{u_-}{2}$ in \eqref{2.3} or \eqref{2.4}, that is, $f^\prime(u_-)>\sigma=f^\prime(u_+)$,  Oleinik shock wave is degenerate at $u_+$,  and we call it degenerate (or sonic) Oleinik shock, and the other case of \eqref{2.3} or \eqref{2.4} is called non-degenerate (or subsonic) Oleinik shock wave. For definiteness, we only consider the case \eqref{2.3}, and the case \eqref{2.4} can be treated similarly.
  
 Correspondingly, the viscous Oleinik shock  is the traveling wave $u^S(\x=x-\sigma t;u_-,u_+)$ defined by the following ODE:
  \begin{equation}\label{2.7}
\left\{\begin{array}{l}
-\sigma(u^S)^\prime+(f(u^S))^\prime=\mu (u^S)^{\prime\prime},\qquad^\prime=\frac{d}{d\xi}, \\[2mm]
\displaystyle u^S(\pm\infty)=u_\pm.
\end{array}\right.
\end{equation}
Now we only consider the degenerate viscous shock case, that is,  $u_+=-\frac{u_-}{2}$ in \eqref{2.3}.  Integrating $ \eqref{2.7}_1 $ from $(-\infty,\xi]$, it holds that
\begin{equation}\label{2.8}
\left\{\begin{array}{l}
\di \mu u^S_\xi=(u^S-u_-)(u^S-u_+)^2,\\[2mm]
\displaystyle u^S(\pm\infty)=u_\pm.
\end{array}\right.
\end{equation}
The existence of the viscous Oleinik shock wave $u^S(\xi)$ can be found in \cite{MK}. The properties of  degenerate viscous Oleinik shock can be summarized in the following lemma, whose proof can be found in \cite{M}.
\begin{lemma}\label{L2-1} 
For any states $u_\pm$ satisfying $u_- <0$ and $u_+=-\frac{u_-}{2}>0$ in \eqref{2.3} and let $\delta_S:=|u_+-u_-|$ be the wave strength of the degenerate Oleinik shock. Then there exists a constant $C>0$ such that the following properties hold true: 
 \begin{equation*}
 u^S_\xi>0,\qquad  \forall \xi\in \mathbb{R},
 \end{equation*}
 and 
\begin{equation}
\begin{aligned}
&\di |u^S(\xi)-u_-|\leq C \delta_S e^{-C\delta_S^2|\xi|},&\di \mbox{if}\ \xi<0,  \notag \\[2mm]
&\di |u^S(\xi)-u_+|\leq\frac{C\delta_S}{1+C\delta_S^2|\xi|},&\di \mbox{if}\  \xi>0, \notag \\[2mm]
&\di |u^S_\xi| \leq C\delta_S^3e^{-C\delta_S^2|\xi|},&\di \mbox{if}\ \xi<0, \notag \\[2mm]
&\di |u^S_\xi| \leq \frac{C \delta_S^3 }{(1+C\delta_S^2|\xi|)^2},&\di \mbox{if}\ \xi>0,  \notag\\[2mm]
&|u^S_{\xi \xi}|\leq C\delta_S^2|u^S_\xi|,&\di \forall \xi\in \mathbb{R}.\notag
\end{aligned}
\end{equation}
\end{lemma}




\

\noindent{\bf $a$-contraction method for a single degenerate Oleinik shock.}   To explain the main ideas of $a$-contraction method applied to the non-convex conservation laws, for simplicity, we fix $u_-=-2$ and $u_+=1$ now and the following arguments still hold true for arbitrary $u_\pm$ satisfying admissible conditions \eqref{2.1} and \eqref{2.2}. Then its degenerate viscous shock speed is given by $\sigma=f^\prime(u_+)=3u_+^2=3$ and the viscous shock profile $U(\xi)$ satisfies
$$
- 3U^\prime+3U^2U^\prime=\mu U^{\prime\prime}, \qquad \lim_{\xi\rightarrow-\infty} U(\xi)=-2, \qquad \lim_{\xi\to+\infty} U(\xi)=1,
$$
where $\xi=x-3t$.
By \eqref{2.8}, we have  
$$
\mu U^\prime (\x)=(U+2)(U-1)^2>0.
$$
Note that the viscous shock profile $U(\xi)$ itself is strictly monotone,  but the derivative of the characteristic field $\partial_{\xi}(f^\prime(U(\xi)))=6U(\xi)U^\prime(\xi)$  can change the sign at the inflection point $U=0$, that is, the strict compressibility of the Oleinik shock profile does not always hold true.

In order to use $a$-contraction method, the following Poincar${\acute{\rm e}}$  type inequality is crucial in our stability proof.

\begin{lemma}\label{L4.1} (\cite{VKM}) For any $f: [0,1]\rightarrow\mathbb{R} $ satisfying $\di \int_0^1 y(1-y)|f^\prime(y)|^2dy<+\infty,$ it holds that
   \begin{equation}\label{4.1}
\int_0^1 \Big|f-\int_0^1 fdy\Big|^2dy\leq\frac{1}{2}\int_0^1y(1-y)|f^\prime(y)|^2dy,
   \end{equation}
   that is,
     \begin{equation}\label{4.1+}
\int_0^1 f^2(y)dy-\Big(\int_0^1 f(y)dy\Big)^2\leq\frac{1}{2}\int_0^1y(1-y)|f^\prime(y)|^2dy.
   \end{equation} 
   \end{lemma}

The shifted Oleinik shock $U^{\mb{X}}(t,\x):=U(\x+{\mathbf{X}}(t))$ satisfies the equation 
\begin{equation*}
\partial_t U^{\mb{X}}-\dot{\mathbf{X}}(t)U^{\mb{X}}_\x- 3U^{\mb{X}}_\x+3(U^{\mb{X}})^2U^{\mb{X}}_\x=\mu U^{\mb{X}}_{\x\x},
\end{equation*}
where the shift function ${\mathbf{X}}(t)$ is defined in \eqref{sf}.
Changing variable $(t,x)\rightarrow (t,\x=x-3t)$, then $u(t,\x):=u(t,\x+3 t)$ satisfies 
\begin{equation*}
u_t-3u_\x+3u^2u_\x=\mu u_{\x\x}.
\end{equation*}
Then the perturbation $\phi(t,\x):=u(t,\x)-U(\x+\mb{X}(t))$ satisfies
\begin{equation}\label{b1}
\phi_t-3\phi_\x+\dot{\mathbf{X}}(t)U^{\mb{X}}_\x+[(\phi+U^{\mb{X}})^3-(U^{\mb{X}})^3]_\x=\mu \phi_{\x\x}.
\end{equation}
Due to the degeneracy and the partial compressibility of viscous Oleinik shock and the influence of the inflection point of non-convex flux, the $a-$contraction and  $L^2$ relative entropy method for Burgers equation in \cite{KV} can not be utilized directly.  We need to use suitable weight in the compressibility region of the shock profile, that is, the negative part of the profile for $U(\xi)<0$.

We denote that  $\xi_1:=(U)^{-1}(0)$ and $\xi_*:=(U)^{-1}(\frac{1}{2})$,  where $(U)^{-1}$ is inverse function of $U$. The weight function $w(U)$ is defined as 
\begin{equation}\label{c8}
w(U(\xi)):=
\begin {cases}
\di \frac{5}{2}(1-U),&\di \quad \xi\in(-\infty,\xi_1)  \iff U\in(-2,0),\\[4mm]
\di \frac{5}{2}(1-U)(4U^3+1),&\di \quad \xi\in[\xi_1,\xi_*)\quad  \iff U\in[0,\frac{1}{2}),\\[4mm]
\di \frac{15}{8},&\di \quad \xi\in[\xi_*,+\infty) \iff U\in[\frac{1}{2},1).\\
\end{cases}
\end{equation}
It is direct to check that the weight function $w(U)$ is $C^2$-smooth and satisfies 
$$\frac{15}{8}\leq w(U)\leq  \frac{15}2.$$
Note that the weight is added essentially on the negative part of the Oleinik shock, that is, on the region $\xi\in(-\infty,\xi_1)$, or equivalently when $U\in(-2,0).$  In the region $[\xi_*,+\infty),$ the weight $w\equiv \frac{15}{8}$ and in fact, no need for weight,  and in the part $[\xi_1,\xi_*)$,  the weight $w$ is chosen such that it is  $C^2$-smooth. We need to point out that this weight $w$ depends on the wave strength such that arbitrarily large shock wave can be handled.

Choose the  time-dependent shift ${\mathbf{X}}(t)$ as
\begin{equation}\label{sf}
\left\{\begin{array}{l}
\di \dot{\mathbf{X}}(t)=\frac{32}{25}\int_{\mathbb{R}} \phi(t,\x) w^{\mb{X}} U^{\mb{X}}_\x d\xi=\frac{32}{25}\int_{\mathbb{R}} \phi^{-\mb{X}}(t,\x) w U_\x d\xi,\\[4mm]
\mathbf{X}(0)=0,
\end{array}\right. 
\end{equation} 
where $w^{\mb{X}}:=w(U^{\mb{X}}(\xi))$ and $\phi^{-\mb{X}}(t,\x) :=\phi(t,\x-\mathbf{X}(t))$.
If we do $L^2$ relative entropy estimate to \eqref{b1} directly, it holds that
\begin{equation}\label{de}
\frac{1}{2}\frac{d}{dt}\int_{\mathbb{R}} \phi^2 d\x+\int_{\mathbb{R}} \phi_\x^2 d\x  +\int_{\mathbb{R}}3U^{\mb{X}}U^{\mb{X}}_\x\phi^2d\x+(\cdots)=0.
\end{equation}
where and in the sequel $(\cdots)$ represents the higher order terms  which can be controlled. Since $U^{\mb{X}}_\x>0,$ it can be seen from \eqref{de} that in the positive region of the Oleinik shock profile $U^{\mb{X}}$, that is,  when $\xi\in[\xi_*,+\infty),$ $L^2$ relative entropy estimate in \eqref{de} is good and no need for the weight and Poincar${\acute{\rm e}}$ inequality, while in the negative region of the Oleinik shock profile $U^{\mb{X}}$, that is,  when $\xi\in(-\infty, \xi_1],$ $L^2$ relative entropy estimate in \eqref{de} has bad sign term and therefore, $a$-contraction method should be used with suitable weight $w$ in \eqref{c8} and time-dependent shift  ${\mathbf{X}}(t)$  in \eqref{sf}.

Precisely, multiplying \eqref{b1} by $w^{\mb{X}}\phi$ with the weight function  $w^{\mb{X}}:=w(U^{\mb{X}}(\xi))$ defined in \eqref{c8}, we can obtain
\begin{equation}
\begin{aligned}
\frac{1}{2}\frac{d}{dt}\int_{\mathbb{R}} \phi^2 w^{\mb{X}}d\x{+\int_{\mathbb{R}}\Big( {3(w^\prime)^{\mb{X}}}-3(U^{\mb{X}})^2(w^\prime)^{\mb{X}}+3U^{\mb{X}}w^{\mb{X}}-\frac{(w^{\prime\prime})^{\mb{X}}}{2}\mu U^{\mb{X}}_\x\Big)U^{\mb{X}}_\x\phi^2d\x}\\
{+\int_{\mathbb{R}}\mu \phi_\x^2w^{\mb{X}} d\x}{+\dot{\mb{X}}(t)\int_{\mathbb{R}}w^{\mb{X}} U^{\mb{X}}_\x \phi d\x}+{(\cdots)}=0.
\end{aligned}
\end{equation}

We claim that if the initial perturbation $\|u_0-U\|_{H^1(\RR)}$ is suitably small, then it holds that
\begin{equation}\label{c5}
\begin{array}{ll}
\di \mathbf{G}(t):=
{\int_{\mathbb{R}} \mu\phi_\x^2w^{\mb{X}} d\x}
{+\dot{\mb{X}}(t)\int_{\mathbb{R}}w^{\mb{X}} U^{\mb{X}}_\x \phi d\x}\\[5mm]
\di\qquad  +{\int_{\mathbb{R}}\Big( {3(w^\prime)^{\mb{X}}}-3(U^{\mb{X}})^2(w^\prime)^{\mb{X}}+3U^{\mb{X}}w^{\mb{X}}-\frac{(w^{\prime\prime})^{\mb{X}}}{2}\mu U^{\mb{X}}_\x\Big)U^{\mb{X}}_\x\phi^2d\x}\\[5mm]
\di\quad  \geq \frac{5}{16}\int_{\mathbb{R}} \mu \phi_\x^2d\x+\frac{4}{5} \int_{\mathbb{R}} \phi^2 U^{\mb{X}}_\x d\x+\frac{25}{64}|\dot{\mb{X}}(t)|^2.
\end{array}
\end{equation}
and  the following weighted $L^2$-contraction holds true,
$$
\frac{d}{dt}\int_{\mathbb{R}} \phi^2 w^{\mb{X}}d\x\leq 0.
$$

If we first change the variable $\x\rightarrow \x-{\mb{X}}(t)$ in \eqref{c5}, then according to the sign of the Oleinik shock profile $U(\x)$, we can split $\mathbf{G}(t)$ as follows
\begin{equation}\label{c7}
\mathbf{G}(t)=\mathbf{G}_1(t)+\mathbf{G}_2(t)
\end{equation}
where 
\begin{equation}\label{c6}
\begin{aligned}
\mathbf{G}_1(t):=&\di \int_{-\infty}^{\xi_*}\mu w (\phi^{-\mb{X}}_\xi)^2d\xi+\int_{-\infty}^{\xi_*}(\phi^{-\mb{X}})^2 U_\xi \Big(3 w^\prime-3U^2w^\prime+3Uw-\frac{w^{\prime\prime}}{2}\mu U_\xi\Big)d\xi\\
\qquad&\di  +\frac{8}{25}\Big(\int_{-\infty}^{\xi_*} \phi^{-\mb{X}} w U_\xi d\xi\Big)^2,
\end{aligned}
\end{equation}
and
\begin{equation}
\begin{aligned}
\mathbf{G}_2(t):=&\di \int_{\xi_*}^{+\infty}\mu w (\phi^{-\mb{X}}_\xi)^2d\xi+3\int_{\xi_*}^{+\infty}(\phi^{-\mb{X}})^2 U w U_\xi d\xi\\
&\di +\mathbf{ \dot  X}(t)\int_{\mathbb{R}} \phi^{-\mb{X}} w U_\xi d\xi-\frac{8}{25}\Big(\int_{-\infty}^{\xi_*} \phi^{-\mb{X}} w U_\xi d\xi\Big)^2.
\end{aligned}
\end{equation}
The key point is to handle the term $\mathbf{G}_1(t)$. 
We can use the modified $a$-contraction method with suitably weighted Poincar${\acute{\rm e}}$ inequality to obtain
\begin{equation}\label{c2}
\begin{aligned}
\mathbf{G}_1(t)\geq  \frac{5}{16}\int_{-\infty}^{\xi_*} \mu(\phi^{-\mb{X}}_\x)^2d\x+\frac{4}{5} \int_{-\infty}^{\xi_*} (\phi^{-\mb{X}})^2 U_\x d\x.
\end{aligned}
\end{equation}
To prove \eqref{c2},  we first let 
\begin{equation}\label{transform}
y:=\frac{U(\x)-u_-}{U(\x_*)-u_-}=\frac{2}{5}(U(\x)+2),
\end{equation}
then we have $\x\in (-\infty,\x_*]\iff y\in[0,1]$, which means that we only use weighted Poincar${\acute{\rm e}}$ inequality on the part $\x\in (-\infty,\x_*]$, not on the whole line $\xi\in \mathbb{R}$,  the main difference from the $a$-contraction method to the classical Burgers equation in \cite{KV}. Since $\di y_\x=\frac{2}{5}U^\prime(\x)>0$, we have from the inverse function theorem and \eqref{transform} that there exists a unique inverse function $\x=\xi(y)$. For any fixed $t>0$, we can write $\psi(t, y): =\phi^{-\mb{X}}(t,\x) w(U(\x))$. By using the weighted Poincar${\acute{\rm e}}$ inequality in Lemma \ref{L4.1} with $\psi(t, y)$, it holds that
\begin{equation}\label{c3}
\begin{aligned}
&\frac{8}{25}\Big(\int_{-\infty}^{\xi_*} \phi^{-\mb{X}} w U_\xi d\xi\Big)^2=2\Big(\int_0^1 \psi(t,y)dy\Big)^2\\
&\geq 2\int_0^1 \psi^2(t,y)dy-\int_0^1y(1-y)|\psi_y(t,y)|^2dy\\
&=-\frac{2}{5}\int_{-\infty}^{\xi_*}\mu (\phi^{-\mb{X}}_\x)^2w^2\frac{\frac{1}{2}-U}{(1-U)^2}d\x
\\&\quad +\frac{2}{5}\int_{-\infty}^{\xi_*} (\phi^{-\mb{X}})^2U_\x\left[ww^{\prime\prime}(\frac{1}{2}-U)(U+2)+2w^2-ww^\prime(2U+\frac{3}{2})\right]d\x.
\end{aligned}
\end{equation}
One can refer to Lemma  \ref{L4.3} for more details. By \eqref{c6} and \eqref{c3}, the estimation of  $\mathbf{G}_1(t)$ in \eqref{c2} can be proved.
Note that the choosing of $\xi_*:=(U)^{-1}(\frac{1}{2})$ is not essential, but crucially we need to choose $\x_*$ such that $U(\x_*)>0$ which means that $U(\x_*)$ should cross the inflection point $0$ of the cubic non-convex flux. Since  the term $\di 3\int_{\xi_*}^{+\infty}(\phi^{-\mb{X}})^2 U w U_\xi d\xi$ has good sign, with the help of the time-dependent shift ${\mathbf{X}}(t)$ defined  in \eqref{sf}, we can get
\begin{equation}\label{c4}
\begin{aligned}
\mathbf{G}_2(t)\geq  \frac{15}{8}\int_{\xi_*}^{+\infty}\mu  (\phi_\xi^{-\mb{X}})^2d\xi+(\frac{45}{16}-\frac{9}{8}\ln 2)\int_{\xi_*}^{+\infty}(\phi^{-\mb{X}})^2 U_\xi d\xi+\frac{25}{64}|\mathbf{ \dot  X}(t)|^2.
\end{aligned}
\end{equation}
By \eqref{c2} and \eqref{c4}, we can prove \eqref{c5}, and then weighted $L^2$-contraction holds 
$$\frac{d}{dt}\int_{\mathbb{R}} \phi^2 w^{\mb{X}}d\x\leq 0.$$

\noindent{\bf  Rarefaction wave.}  We first recall  the Riemann problem for the inviscid Burgers equation:
  \begin {equation}\label{2.12}
  \left\{\begin{array}{l}
  w^r_t+w^rw^r_x=0, \\
 w^r(0,x)=w^r_0(x)=  \left\{\begin{array}{l}
 w_-, x<0,\\
 w_+,x>0.
 \end{array}\right.\\
  \end{array}\right.\\
   \end {equation}
If $w_-<w_+$, then the Riemann problem \eqref{2.12} has a self-similar rarefaction wave fan solution $\di w^r(t,x):=w^r(\frac{x}{t};w_-,w_+)$ defined by
\begin{equation}\label{2.13}
w^r(t,x)=w^r(\frac{x}{t};w_-,w_+) :=
\left\{\begin{array}{l}
w_-,\ \ \ x\leq w_-t,\\[2mm]
\di \frac{x}{t},\ \ \ w_-t\leq x\leq w_+t,\\[2mm]
w_+,\ \ \ x\geq w_+t.\\
\end{array}\right. \\
\end{equation}
For the Riemann problem of general scalar hyperbolic conservation law \eqref{1.2}, if $f^{\prime\prime}(u)>0$ on $u\in[u_-,u_+]$ and $f^\prime(u_-)<f^\prime(u_+)$, then the self-similar rarefaction wave solution $\di u^r(t, x):=u^r\big(\frac{x}{t}; u_-,u_+\big)$ can be defined explicitly by
\begin{equation}\label{2.14}
u^r\big(\frac{x}{t}; u_-,u_+\big)=(\lambda)^{-1}\Big(w^r\big(\frac{x}{t};\lambda_-,\lambda_+\big)\Big),
\end{equation}
where $\lambda(u):=f^\prime (u)$ and $\lambda_\pm:=\lambda(u_\pm)=f^\prime(u_\pm)$.

\noindent{\bf Composite wave of degenerate Oleinik shock and rarefaction.} Given the Riemann initial data $u_\pm$ , we consider the cubic Burgers equation, i.e. flux $f(u)=u^3$. If for $u_-<0,$ the state $u_+$ satisfies that
    \begin{equation}\label{2.17}
  u_-<-\frac{u_-}{2}<u_+,
  \end{equation}
  or for $u_->0,$ 
    \begin{equation}\label{2.18}
  u_+<-\frac{u_-}{2}<u_-,
  \end{equation}
  then  the Riemann problem \eqref{1.2} is solved by composite wave of  a degenerate Oleinik shock connecting $u_-$ and $u_m(=-\frac{u_-}{2})$ and a rarefaction wave connecting $u_m$ and $u_+$. Correspondingly,  the time-asymptotic stability of the Cauchy problem \eqref{1.1} is conjectured to be determined by composite wave of a degenerate viscous Oleinik shock and a rarefaction wave. 
  
  Without loss of generality, in  the sequel we consider the far fied states $u_\pm$ satisfy \eqref{2.17}, and we always denote $u_m:=-\frac{u_-}{2}$. Then the large-time behavior of the Cauchy problem \eqref{1.1}  satisfying \eqref{2.17}  is determined by  a superposition wave:
  \begin{equation}\label{2.19}
 \Big(u^S(x-\sigma t+ \mathbf{X}(t); u_-,u_m)+u^r(\frac{x}{t};u_m,u_+)-u_m\Big),
  \end{equation}
  where  $u^S(x-\sigma t+ \mathbf{X}(t); u_-,u_m)$ is defined in \eqref{2.7}  and  $\di u^r(\frac{x}{t};u_m,u_+)$ is defined in \eqref{2.14}.
  
\noindent{ \bf Main result.}  Now we can state our main result as follows.
\begin{theorem}\label{T2.1} If the far-fields states $u_{\pm}$ in \eqref{1.1} satisfy \eqref{2.17}, then there exist positive constants  $\delta^*$, $\epsilon^*$ such that if the initial data $u_0$ satisfies 
 \begin{equation}\label{2.20}
 \begin{array}{ll}
\di \big{ \|}u_0(\cdot)-u^S(\cdot;u_-,u_m)\big{\|}_{L^2(\mathbb{R_-})}+ \big{ \|}u_0(\cdot)-\big(u^S(\cdot;u_-,u_m)+u_+-u_m\big)\big{\|}_{L^2(\mathbb{R_+})}\\[3mm]
\hspace{7.32cm}\di+\big{ \|}u_{0x}(\cdot)-u^S_x(\cdot;u_-,u_m)\big{\|}_{L^2(\mathbb{R})}<\epsilon^*,
\end{array}
  \end{equation}
 with $\mathbb{R_+}:=(0,+\infty)=:-\mathbb{R_-}$, and the rarefaction wave strength  $\delta_R$ satisfies $\delta_R\leq\delta^*$, 
then the Cauchy problem \eqref{1.1} admits a unique global-in-time solution $u(t,x)$ satisfying 
 \begin{equation}\label{2.21}
 \begin{aligned}
&u(t,x)- \Big(u^S\big(x-\sigma t + \mathbf{X}(t);u_-,u_m\big)+u^r\Big(\frac{x}{t};u_m,u_+\Big)-u_m\Big)\in C(0,+\infty; H^1(\mathbb R)),\\
&u_{xx}(t,x)- u^S_{xx}\big(x-\sigma t + \mathbf{X}(t);u_-,u_m\big) \in L^2(0,+\infty; L^2(\mathbb R)),
\end{aligned}
\end{equation}
where the shift function $\mathbf{X}(t)$ is absolutely continuous. Moreover, the following time-asymptotic stability of composite wave of degenerate Oleinik shock with shift  $\mathbf{X}(t)$ and rarefaction wave holds true,
 \begin{equation}\label{2.22}
 \displaystyle\lim_{t\rightarrow+\infty} \displaystyle \sup_{x\in\mathbb R}\Big|u(t,x)- \Big(u^S\big(x-\sigma t + \mathbf{X}(t);u_-,u_m)+u^r(\frac{x}{t};u_m,u_+)-u_m\Big)\Big|= 0,
 \end{equation}
 and the shift function $\mathbf{X}(t)$ satisfies
 \begin{equation}\label{2.23}
 \displaystyle\lim_{t\rightarrow+\infty}|\dot{\mathbf{X}}(t)|= 0.
 \end{equation}
 \end{theorem} 
 
 \
 
\noindent {\bf Remark 2.1}. The shift ${\mathbf{X}}(t)$ is proved to satisfy the time-asymptotic behavior \eqref{2.23}, which implies 
   \begin{equation}\label{2.24}
 \displaystyle\lim_{t\rightarrow+\infty}\frac{\mathbf{X}(t)}{t}= 0, \notag
 \end{equation}
 that is, the shift $\mathbf{X}(t)$ grows at most sub-linearly with respect to the time t, therefore, the shifted viscous Oleinik shock $u^S\big(x-\sigma t + \mathbf{X}(t);u_-,u_m)$ still keeps the original Oleinik shock profile time-asymptotically.
 
 \
  
\noindent{\bf Remark 2.2}.  Theorem \ref{T2.1} is the first time-asymptotic stability result towards the  composite wave of Oleinik shock and rarefaction wave to scalar non-convex conservation laws as far as we know, which is proposed as an open problem in the survey paper of Matsumura \cite{M}.
Moreover, in Theorem \ref{T2.1}, Oleinik shock wave strength can be arbitrarily large, while the rarefaction wave strength should be suitably small due to the subtle wave interactions.

 \
  
\noindent{\bf Remark 2.3}.   Note that both the initial value $u_0$ and the initial Oleinik shock $u^S$ in the initial perturbation \eqref{2.20} belong to $H^1(\mathbb{R})$ up to the constant spatial far-fields, while the initial self-similar rarefaction wave is the discontinuous Riemann data and belongs to  piece-wise $H^1(\mathbb{R_{\pm}})$.



\section{Reformulation of the problem}   
In this section, we first construct an approximate rarefaction wave. Then we reformulate the stability problem around the shifted viscous Oleinik shock and shifted approximate rarefaction wave based on the construction of the time-dependent shift $\mb{X}(t)$ and suitable weight function  $w(u^S(\x))$.
Finally, we give the proof of Theorem \ref{T2.1} based on the local existence of the solution and the uniform-in-time a priori estimates.

\subsection{Construction of approximate rarefaction wave}

 Note that the self-similar rarefaction wave fan solutions $\di w^r(t,x)$ and $\di u^r(t,x)$ in \eqref{2.13} and \eqref{2.14}
are Lipschitz continuous for $t>0$ and $x\in \RR$. In order to study their time-asymptotic stability for the viscous equation \eqref{1.1} with second order derivative, we need to construct their smooth approximations which are equivalent to the  corresponding  inviscid self-similar rarefaction waves time-asymptotically as in \cite{MN1}. Precisely, consider the following Cauchy problem of Burgers equation:
 \begin {equation}\label{2.15}
  \left\{\begin{array}{l}
\di  w^R_t+w^Rw^R_x=0, \\[2mm]
\di  w^R(0,x)=w^R_0(x)= \frac{w_++w_-}{2}+\frac{w_+-w_-}{2}\tanh x\ra w_{\pm},\ \mbox{as}\ \ x\ra\pm\infty.
  \end{array}\right.\\
   \end {equation}
 Since $w_-<w_+$ and $(w^R_0)^\prime(x)>0$ for any $x\in\RR$, the Cauchy problem \eqref{2.15} has a unique smooth solution $\di w^R(t,x):=w^R(t,x;w_-,w_+)$ defined by 
\begin{equation}\label{a2}
\left\{\begin{array}{l}
\di w^R(t,x)=w^R_0\big(x_0(t,x)\big)=\frac{w_++w_-}{2}+\frac{w_+-w_-}{2}\tanh\big(x_0(t,x)\big),\\[3mm]
\di x=x_0(t,x)+w^R_0\big(x_0(t,x)\big)t,
  \end{array}\right.
    \end{equation}
where $x_0(t,x)$ is the unique intersection point with $x$-axis of the straight characteristic line of Burgers equation passing through the point $(t,x)$ for $t>0$ and $x\in \RR$.
Correspondingly, a smooth approximate rarefaction wave $u^R(t,x):=u^R(t,x; u_-,u_+)$ of self-similar rarefaction wave $\di u^r(t, x):=u^r\big(\frac{x}{t}; u_-,u_+\big)$ can be constructed by
 \begin {equation}\label{2.16}
 u^R(t,x; u_-,u_+):=(\lambda)^{-1}\big(w^R(t,x;\lambda_-,\lambda_+)\big),
  \end {equation}
with $\lambda(u):=f^\prime (u)$ and $\lambda_\pm:=\lambda(u_\pm)=f^\prime(u_\pm)$. Then $u^R(t,x)$ in \eqref{2.16} satisfies the equation
 \begin{equation}\label{arare}
\left\{\begin{array}{l}
u_t^R+f(u^R)_x=0,\\[2mm]
u^R(0,x)=(\lambda)^{-1} \Big(\frac{\lambda_++\lambda_-}{2}+\frac{\lambda_+-\lambda_-}{2}\tanh x\Big)\ra u_{\pm},\ \mbox{as}\ \ x\ra\pm\infty.
\end{array}\right. 
\end{equation}
Moreover, we have the following lemma.

\begin{lemma}\label{L2-2} 
Assume $u_-<u_+$ and $f^{\prime\prime}(u)>0$ on $u\in [u_-,u_+]$ and let $\delta_R:=|u_+-u_-|$ be the rarefaction wave strength. Then  it holds that 

(1). $ u_-<u^R(t,x)<u_+$, $ u_x^R(t,x)>0$, for $t\geq0$ and  $x\in\mathbb{R}$;

(2). For all $t\geq1$ and $p\in[1,+\infty]$, there exists a positive constant $C_p$ such that
\begin{equation*}
 \begin{array}{ll}
\di  \|u_{x}^R(t)\|_{L^p(\mathbb{R})}\leq C_p \min\left\{ \delta_R,{\delta_R}^{\frac{1}{p}}t^{-1+\frac{1}{p}}\right\},\\[3mm]
\di \|u_{xx}^R(t)\|_{L^p(\mathbb{R})}\leq C_p \min\left\{\delta_R,t^{-1}\right\};
\end{array}
\end{equation*}
 
 (3). $\di \lim_{t\rightarrow \infty} \sup_{x\in\mathbb{R}}|u^R(t,x)-u^r(\frac{x}{t})|=0$, that is, the approximate rarefaction wave $u^R(t,x)$ and the inviscid self-similar rarefaction wave fan $\di u^r(\frac xt)$ are equivalent to each other time-asymptotically for $x\in \RR$ uniformly;
 
 (4). For all  $t\geq0$, it holds that
  \begin{equation}
|u^R(t,x)-u_+|\leq C\delta_{R}e^{-2|x-\lambda_+t|},\ \ x\geq\lambda_+t,\notag
\end{equation}
  \begin{equation}
|u^R(t,x)-u_-|\leq C\delta_{R}e^{-2|x-\lambda_-t|},\ \ x\leq\lambda_-t;\notag
\end{equation}

 (5). For all  $t\geq1$ and  any $\ve\in(0,1)$, there exists a positive constant $C_\ve$ such that 
 \begin{equation}
 |u^R(t,x)-u_+|\leq C_{\ve}\delta_{R}^{\frac{2\ve}{2+\ve}}t^{-1+\ve}e^{-\ve|x-\lambda_+t|}, \ \ x\geq\lambda_+t,\notag
\end{equation}
 \begin{equation}
 |u^R(t,x)-u_-|\leq C_{\ve}\delta_{R}^{\frac{2\ve}{2+\ve}}t^{-1+\ve}e^{-\ve|x-\lambda_-t|}, \ \ \ x\leq\lambda_-t; \notag
 \end{equation}
 
 (6). For all  $t\geq1$ and  any $\ve\in(0,1)$, there exists a positive constant $C_\ve$ such that 
 \begin{equation}
 |u^R(t,x)-u^r(\frac{x}{t})|\leq C_{\ve}\delta_{R}^{\ve} t^{-1+\ve},\ \  \lambda_-t\leq x\leq\lambda_+t. \notag
\end{equation}

\end{lemma}

\begin{proof} The proof of properties (1)-(3) can be found exactly in \cite{MN1}, and (4) can be proved directly by the solution formula \eqref{a2} and \eqref{2.16}. Now we prove (5) and (6), whose proofs are motivated by \cite{MY}, but here we need more detailed information about the rarefaction wave strength $\delta_R$. For simplicity, 
we only prove (6) and the second inequality in (5), while the first one in (5) can be done similarly. 

By \eqref{2.16}, to prove the second inequality in (5), it is enough to show
  \begin{equation}\label{a1}
  |w^R(t,x)-\lambda_-|\leq C_{\ve}\delta_{R}^{\frac{2\ve}{2+\ve}}t^{-1+\ve}e^{-\ve|x-\lambda_-t|},\ \  \forall t\geq1\ \mbox{and}\  x\leq\lambda_-t,
    \end{equation}
for the solution $w^R(t,x)$ in \eqref{2.15} with $w_{\pm}=\lambda_{\pm}:=f^\prime(u_\pm)$.
It follows from \eqref{a2} that
    \begin{equation}\label{a3}
\frac{\partial x_0}{\partial x}=\frac{1}{1+(w^R_0)^\prime(x_0)t}>0,
    \end{equation}
which implies that
  \begin{equation}\label{a4}
x_0(t,x)\leq x_0(t,\lambda_-t), \ \ \forall x\leq\lambda_-t.
    \end{equation}
Again by \eqref{a2},  $x_0(t,\lambda_-t)$ is given by
   \begin{equation}\label{a5}
 x_0(t,\lambda_-t)=\lambda_-t-w^R_0 \big(x_0(t,\lambda_-t)\big)t<0, \ \ \forall t>0,
    \end{equation}
which results in 
 \begin{equation}\label{a6}
\frac{d}{dt} x_0(t,\lambda_-t)=\frac{\lambda_--w^R_0 \big(x_0(t,\lambda_-t)\big)}{1+(w_0^R)^\prime \big(x_0(t,\lambda_-t)\big) t}<0,\ \ \forall t>0.
    \end{equation}
Therefore, by \eqref{a5} and  \eqref{a6}, as $t\rightarrow +\infty,$
 \begin{equation}\label{a7}
 x_0(t,\lambda_-t)\rightarrow -\infty.
    \end{equation}
  By \eqref{a2} and \eqref{a5}, we have
   \begin{equation}\label{a8}
 x_0(t,\lambda_-t)=(\lambda_--\lambda_+)\frac{e^{ 2x_0(t,\lambda_-t)}}{e^{ 2x_0(t,\lambda_-t)}+1}t<0, \ \ \forall t>0.
    \end{equation}
  Thus,  there exist uniform constants $C_1>C>0$  such that 
 \begin{equation}\label{a9}
-C_1 \delta_R\ t\ e^{ 2x_0(t,\lambda_-t)}\leq x_0(t,\lambda_-t)\leq -C\delta_R\ t\ e^{ 2x_0(t,\lambda_-t)}, \ \ \forall t>0.
    \end{equation}  
Therefore, for any $\ve>0$, based on the fact $x\leq C_\ve e^{\ve x}$, we have
 \begin{equation}\label{a10}
\delta_R\ t\leq C_{\ve} e^{-(2+\ve) x_0(t,\lambda_-t)},
    \end{equation}  
  which implies
   \begin{equation}\label{a11}
  C_{\ve}(\delta_R t)^{-\frac{2}{2+\ve}}\geq e^{2x_0(t,\lambda_-t)}.
   \end{equation}   
   Hence by \eqref{a2}, \eqref{a5}, \eqref{a9} and \eqref{a11}, we have $\forall  x\leq \lambda_-t,$
     \begin{equation}\label{a12}
 \begin{aligned}
|w^R(t,x)-\lambda_-|=w^R_0\big(x_0(t,x)\big)-\lambda_-
\leq w^R_0\big(x_0(t,\lambda_-t)\big)-\lambda_-\\
=-\frac{x_0(t,\lambda_-t)}{t}\leq  C\delta_R e^{ 2x_0(t,\lambda_-t)}
\leq   C_{\ve} \delta_R^{\frac{\ve}{2+\ve}}t^{-\frac{2}{2+\ve}}.
   \end{aligned}
    \end{equation}
On the other hand, since
 \begin{equation}\label{a13}
x=x_0(t,x)+w^R_0(x_0(t,x))t\geq x_0(t,x)+\lambda_-t,
   \end{equation}   
    we have for $x\leq \lambda_-t,$
    \begin{equation}\label{a14}
 \begin{aligned}
 |w^R(t,x)-\lambda_-|=w^R_0(x_0(t,x))-\lambda_-=(\lambda_+-\lambda_-)\frac{e^{2x_0(t,x)}}{e^{2x_0(t,x)}+1}\\[1mm]
 \leq C\delta_R e^{2x_0(t,x)}\leq  C\delta_R e^{2(x-\lambda_-t)}=  C\delta_R e^{-2|x-\lambda_-t|}.
   \end{aligned}
   \end{equation}  
   By \eqref{a12} and  \eqref{a14}, we have for $\forall  t\geq1$ and $\forall x\leq \lambda_-t$,
    \begin{equation}\label{a15}
 \begin{aligned}
 |w^R(t,x)-\lambda_-|= |w^R(t,x)-\lambda_-|^{1-\frac{\ve}{2}} |w^R(t,x)-\lambda_-|^{\frac{\ve}{2}}  \\
 \leq C_\ve\delta_R^{\frac{2\ve}{2+\ve}}t^{-\frac{2-\ve}{2+\ve}}e^{-\ve|x-\lambda_-t|}
  \leq  C_\ve \delta_R^{\frac{2\ve}{2+\ve}}t^{-1+\ve}e^{-\ve|x-\lambda_-t|}.
  \end{aligned}
   \end{equation}  
  Now we prove (6). From \eqref{a2}, we have
  \begin{equation}\label{a16}
  \frac{x}{t}-w^R(t,x)=\frac{x_0(t,x)}{t}.
  \end{equation} 
  Then, by \eqref{a3}, we have for $\lambda_-t\leq x\leq \lambda_+t,$
    \begin{equation}\label{a17}
\frac{x_0(t,\lambda_-t)}{t}\leq  \frac{x}{t}-w^R(t,x)\leq \frac{x_0(t,\lambda_+t)}{t}.  \end{equation} 
Thus, we have for $\lambda_-t\leq x\leq \lambda_+t$,
 \begin{equation}\label{a18}
\big|\frac{x}{t}-w^R(t,x)\big|\leq \max \left\{\Big|\frac{x_0(t,\lambda_-t)}{t}\Big|, \Big|\frac{x_0(t,\lambda_+t)}{t}\Big|\right\}.  \end{equation} 
By \eqref{a9} and  \eqref{a11}, we get 
  \begin{equation}\label{a19}
  \Big|\frac{x_0(t,\lambda_-t)}{t}\Big|\leq C \delta_{R} e^{2x_0(t,\lambda_-t)}\leq C_{\ve}\delta_{R}^{1-\frac{2}{2+\ve}}t^{-\frac{2}{2+\ve}},\ \  \forall \ve>0,\ \forall t\geq 1.
  \end{equation}
Similarly, we have
 \begin{equation}\label{a20}
   \Big|\frac{x_0(t,\lambda_+t)}{t}\Big|\leq C_{\ve}\delta_{R}^{1-\frac{2}{2+\ve}}t^{-\frac{2}{2+\ve}},\  \forall \ve>0, \ \forall t\geq 1.
   \end{equation}
Hence, we have  $\forall \ve\in (0,1), \forall t\geq 1,\  \di  \left|u^R(t,x)-u^r(\frac{x}{t})\right|\leq C_{\ve}\delta_{R}^{\ve} t^{-1+\ve}$ for $\lambda_-t\leq x\leq\lambda_+t$.
 \end{proof}

\subsection{Reformulation of the problem} For the Cauchy problem \eqref{1.1} with far field states $u_\pm$ satisfying \eqref{2.17}, we want to prove its time-asymptotic ansatz is
   \begin{equation}\label{3.1}
   \bar u(t,x):= u^S\big(x-\sigma t + \mathbf{X}(t);u_-,u_m)+u^r\big(\frac{x}{t};u_m,u_+\big)-u_m,
   \end{equation}
   with the time-dependent shift function $\mb{X}(t)$ to be determined. Since the self-similar rarefaction wave $\di u^r\big(\frac{x}{t};u_m,u_+\big)$ is only Lipschitz continuous, in the subsequent stability analysis, the rarefaction wave $\di u^r\big(\frac{x}{t};u_m,u_+\big)$
is replaced by a shifted approximation rarefaction $u^R\big(1+t,x+\mathbf{X}(t);u_m,u_+\big)$ constructed in Section 3.1,  with the exactly same time-dependent shift function $\mb{X}(t)$ as viscous Oleinik shock, such that the stability ansatz is defined by
    \begin{equation}\label{3.2}
    \tilde u(t,x+ \mathbf{X}(t)):=u^S\big(x-\sigma t + \mathbf{X}(t); u_-, u_m\big)+u^R\big(1+t,x+\mathbf{X}(t);u_m,u_+\big)-u_m.
     \end{equation}
For simplification, we denote 
 \begin{equation}\label{3.2+}
    \tilde u^{\mb{X}}(t, x):=\tilde u(t,x+ \mathbf{X}(t)).
     \end{equation}
Note that the above defined ansatz is  quite different from the case of the time-asymptotic stability of the composite waves including viscous shock wave and rarefaction wave with the time-dependent shift function only for the viscous shock wave to the barotropic compressible Navier-Stokes equations in \cite{KVW}, such that the shifted viscous shock wave and the rarefaction wave are still separate in different characteristic fields and can not pass through each other, while here the Oleinik shock and rarefaction wave are always attached together for scalar non-convex conservation laws. Therefore, the same time-dependent shift function need to be equipped to both Oleinik shock and rarefaction wave such that these two waves can not pass through each other and the wave interactions are  still good enough.   

For simplification of our analysis, we rewrite the cubic Burgers equation \eqref{1.1} through the  coordinates transformation $(t,x)\rightarrow(t,\xi=x-\sigma t)$, based on the change of variable associated to the speed of propagation of the shock, $u(t,\xi):=u(t,\xi+\sigma t)=u(t, x)$ satisfies
\begin{equation}\label{1.1+}
u_t-\sigma u_\x+f(u)_\x=\mu u_{\x\x}.
\end{equation}


Meanwhile, in the  coordinates $(t,\xi=x-\sigma t)$, 
 the ansatz $\tilde u(t,x+\mathbf{X}(t))=\tilde u(t, \x+\sigma t+\mathbf{X}(t))=:\tilde u^{\mb{X}}(t,\x)$ in \eqref{3.2} satisfies the equation
\begin{equation}\label{3.6}
\begin{array}{ll}
\di \tilde u^{\mb{X}}_t-\sigma \tilde u^{\mb{X}}_\xi +f(\tilde u^{\mb{X}})_\xi-\mu \tilde u^{\mb{X}}_{\xi\xi}
 =F^{\mb{X}}+ \mathbf{\dot X}(t)  (u^R)^{\mb{X}}_\xi+ \mathbf{\dot X}(t) \ (u^S)^{\mb{X}}_\xi,
\end{array}
\end{equation}
where and in the sequel $\di (u^R)^{\mb{X}}:=u^R(1+t,\x+\sigma t+\mathbf{X}(t)) $, $\di (u^S)^{\mb{X}}:=u^S(\x+\mathbf{X}(t))$ and the error term $F^{\mb{X}}$ is given by
\begin{equation}\label{3.7}
\begin{array}{ll}
\di F^{\mb{X}}=\big[f(\tilde u^{\mb{X}})-f((u^R)^{\mb{X}})-f((u^S)^{\mb{X}})\big]_\x
-\mu (u^R)^{\mb{X}}_{\xi\xi}\\[3mm]
\di\quad = \Big[\big(f^\prime(\tilde u^{\mb{X}})-f^\prime((u^S)^{\mb{X}})\big)(u^S)^{\mb{X}}_\xi+\big(f^\prime(\tilde u^{\mb{X}})-f^\prime((u^R)^{\mb{X}})\big)(u^R)^{\mb{X}}_\xi\Big]-\mu (u^R)^{\mb{X}}_{\xi\xi}\\[3mm]
\di\quad :=F_1^{\mb{X}}+F_2^{\mb{X}},
\end{array}
\end{equation}
with $F_1^{\mb{X}}$ being the wave interactions and $F_2^{\mb{X}}$ the error due to the inviscid rarefaction wave profile. For simplification, we also denote 
 \begin{equation}\label{3.7+}
    F:= \Big[\big(f^\prime(\tilde u)-f^\prime(u^S)\big)u^S_\xi+\big(f^\prime(\tilde u)-f^\prime(u^R)\big)u^R_\xi\Big]-\mu u^R_{\xi\xi},
     \end{equation}
i.e. $ F^{\mb{X}}(t,\xi)=F(t,\xi+\mathbf{X}(t))$.

We set the perturbation
\begin{equation}\label{3.9}
\phi(t,\xi):=u(t,\x)-\tilde u^{\mb{X}}(t, \x).
\end{equation}
Then by \eqref{1.1+} and  \eqref{3.6}, the perturbation $\phi$ satisfies the equation
 \begin{equation}\label{3.10}
\left\{\begin{array}{l}
\phi_t-\sigma\phi_\xi +\big(f(\phi+\tilde u^{\mb{X}})-f(\tilde u^{\mb{X}})\big)_\xi+\mathbf{\dot X}(t)((u^S)^{\mb{X}}_\xi+(u^R)^{\mb{X}}_\xi)-\mu \phi_{\xi\xi}=-F^{\mb{X}},\\[3mm]
\phi(0,\xi)=\phi_0(\xi):=u_0(\xi)-\tilde u(0,\xi),
\end{array}\right.
\end{equation}
with $F^{\mb{X}}$ defined in \eqref{3.7}.

\subsection{Construction of weight function}
   
  By  Lemma \ref{L2-1}: $u^S_\xi>0$, then there exists a unique  $\xi_1$ and a unique $\xi_*$ such that $u^S(\xi_1)=0$ and $u^S(\xi_*) =u_*:=\frac{u_m}{2}$, respectively.  We define the weight function $w=w(u^S(\x))$ by
  \begin{equation}\label{4.2}
w(u^S(\x)):=
\begin {cases}
\di \frac{5}{2}u_m(u_m-u^S(\x)),&\di \quad \xi\in(-\infty,\xi_1)  \iff u^S\in(u_-,0),\\[4mm]
\di \frac{5}{2u_m^2}(u_m-u^S)[4(u^S)^3+u_m^3],&\di \quad \xi\in[\xi_1,\xi_*)\quad  \iff u^S\in[0,u_*),\\[4mm]
\di \frac{15}{8}u_m^2,&\di \quad \xi\in[\xi_*,+\infty) \iff u^S\in[u_*,u_m).\\
\end{cases}
\end{equation}
It is easy to check that $w\in C^2(\mathbb{R})$. 
 Notice that
   \begin{equation}\label{4.3}
   \frac{15}{8}u_m^2\leq w<\frac{15}{2}u_m^2,
   \end{equation}
 \begin{equation}\label{4.4}
 -\frac{5}{2}u_m\leq w^\prime\leq 0\ \ \ \mbox{and}\ \ \  0\leq w^{\prime\prime} \leq \frac{15}{2},\qquad w^{\prime}:=\frac{d w(u^S)}{d u^S}.
   \end{equation}
 This weight function $w$ plays an important role in our stability proof.

\subsection{Construction of shift function}
 Define the shift function $\mathbf{X}(t)$ as a solution to the ODE:
\begin{equation}\label{3.11}
\left\{\begin{array}{l}
\di \dot{\mathbf{X}}(t)=\frac{32}{25u_m^2}\int_{\mathbb{R}} \phi(t,\x) w^{\mb{X}}(u^S(\x)) (u^S)^{\mb{X}}_\xi(\x) d\xi,\\[4mm]
\mathbf{X}(0)=0,
\end{array}\right. 
\end{equation}  
where
$$
(u^S)^{\mb{X}}(\x)=u^S(\x+\mathbf{X}(t)) \ \ \ \mbox{and}\ \ \  w^{\mb{X}}(u^S(\x)):=w((u^S)^{\mb{X}}(\x)).
$$


By applying Lemma A.1 in \cite{KMY} to  \eqref{3.11},  we can prove that \eqref{3.11} has a unique absolutely continuous solution defined on any interval in time $[0,T]$. Moreover, since \eqref{1.1} is uniformly parabolic, the maximum principle \cite{IO} implies that $\displaystyle\sup_{t\in[0,T], x\in\mathbb R} |u(t,x)|\leq C$. By \eqref{3.9}, we have $\displaystyle\sup_{t\in[0,T], \xi \in\mathbb R}|\phi(t,\xi)|\leq C$. Thanks to the facts that $\|w\|_{C^1(\mathbb{R})}\leq C$\ and \ $\|u^S\|_{C^1(\mathbb{R})}\leq C$, we have
 \begin{equation}\label{3.13}
\displaystyle \sup_{\mathbf{X}\in\mathbb{R}}\left|\frac{32}{25u_m^2}\int_{\mathbb{R}} \phi(t,\x) w^{\mb{X}}(u^S(\x)) (u^S)^{\mb{X}}_\xi(\x) d\xi\right|\leq C,
\end{equation}
and
\begin{equation}\label{3.14}
\begin{array}{ll}
\displaystyle \sup_{\mathbf{X}\in\mathbb{R}}\left|\partial_{\mathbf{X}}\left(\frac{32}{25u_m^2}\int_{\mathbb{R}} \phi(t,\x) w^{\mb{X}}(u^S(\x)) (u^S)^{\mb{X}}_\xi(\x) d\xi\right)\right|\\[4mm]
\di \leq C\int_{\mathbb{R}}\big| \big(w^{\mb{X}} (u^S)^{\mb{X}}_\x\big)_\x\big|d \xi \leq C\|w\|_{C^1(\mathbb{R})}\|u^S\|_{C^1(\mathbb{R})}\leq C.
\end{array}
\end{equation}
Especially, since $ |\mathbf{ \dot  X}(t)|\leq C$ by \eqref{3.13}, we have
\begin{equation}\label{3.15}
|\mathbf{X}(t)| \leq Ct, \ \ \  \forall t\leq T.
\end{equation}

\noindent {\bf Proof of Theorem \ref{T2.1}}.  In order to prove Theorem \ref{T2.1}, we shall combine a local existence result together with a priori estimates by continuity arguments.

To state the local existence, we define $\hat U(t,x)=u^S(x-\sigma t;u_-,u_m)+u^R(1+t,x;u_m,u_+)-u_m$.  In the coordinates transformation $(t,x)\rightarrow(t,\x=x-\sigma t)$, $\hat U(t,x)=\hat U(t,\x+\sigma t)=:\hat U(t,\x)$ satisfies
$$
\hat U_t-\sigma\hat U_{\x}+f(\hat U)_{\x}-\mu\hat U_{\x\x}=\hat F,
$$
where  
$$
\hat F=\big(f^\prime(\hat U)-f^\prime(u^S)\big)u^S_\xi+\big(f^\prime(\hat U)-f^\prime(u^R)\big)u^R_\xi-\mu u^R_{\xi\xi}.
$$ 
Set perturbation $\hat \phi(t,\x) =u(t,\x)-\hat U(t,\x)$. By \eqref{1.1+}, we have
$$
\hat \phi_t-\sigma\hat \phi_{\x}+f(\hat \phi+\hat U)_{\x}-f(\hat U)_\x  -\mu\hat \phi_{\x\x}=-\hat F.
$$
Now, we reformulate the Cauchy problem at general initial time $\tau\geq 0$:
 \begin{equation}\label{3.16}
\left\{\begin{array}{l}
\hat \phi_t-\sigma\hat \phi_{\x}+f(\hat \phi+\hat U)_{\x}-f(\hat U)_\x  -\mu\hat \phi_{\x\x}=-\hat F,\\[2mm]
\hat\phi(\tau,\x)=\hat\phi_{\tau}(\x),
\end{array}\right. \ \ \ (t,\x)\in  \mathbb{R}_{+}\times \mathbb{R}.
\end{equation}
\begin{proposition}\label{T3.1}(Local existence) For any $M>0$, there exists a positive constant $t_0=t_0(M)$ independent of $\tau$ such that if $\hat \phi_\tau\in H^1(\mathbb{R})$ and $\|\hat \phi_\tau\|_{H^1(\mathbb{R})}\leq M$, then the Cauchy problem \eqref{3.16} has a unique solution $\hat \phi$ on the time interval $[\tau,\tau+t_0]$ satisfying
\begin{equation}\label{3.17}
\left\{\begin{array}{l}
 \hat \phi\in C([\tau,\tau+t_0];H^1(\RR))\cap L^2(\tau,\tau+t_0; H^2(\RR)),\\[3mm]
 \displaystyle \sup_{t\in[\tau,\tau+t_0]}\| \hat \phi(t)\|_{H^1(\mathbb{R})}\leq2M.
\end{array}\right.
\end{equation}
\end{proposition}
\begin{proposition}\label{T3.2}(A priori estimates) For given $u_-<0$ and the far field states $u_\pm$ satisfying \eqref{2.17}, there exist positive constants $\delta_0,\epsilon_0>0$, such that if the Cauchy problem \eqref{3.10} has a solution $\phi\in C([0,T];H^1(\mathbb{R}))\cap L^2(0,T; H^2(\mathbb{R}))$ for some time $T>0$ with
\begin{equation}
\mathcal{E}(T):=\displaystyle \sup_{0\leq t\leq T}\|\phi\|_{H^1(\mathbb{R})}\leq \epsilon_0, \ \ \  \delta_R:=|u_+-u_m|\leq \delta_0, \notag
\end{equation}
then there exists a uniform-in-time positive constant $C_0$ such that for all $t\in[0,T]$
\begin{equation}\label{3.18}
 \begin{aligned}
\|\phi(t)\|^2_{H^1(\mathbb{R})}+\int_0^t\|\phi_\xi\|^2_{H^1(\mathbb{R})}d\tau+\int_0^t\int_\mathbb{R} {\phi}^2((u^S)^{\mb{X}}_\xi+(u^R)^{\mb{X}}_\xi )d \xi d\tau\\[3mm]
\di+\int_0^t|  \mathbf{ \dot  X}(\tau)|^2d\tau\leq C_0(\|\phi_0\|^2_{H^1(\mathbb{R})}+\delta_R^\frac{8}{33}).
 \end{aligned}
\end{equation}
In addition, by \eqref{3.11},
\begin{equation}\label{3.19}
| \mathbf{ \dot  X}(t)|\leq C\|\phi\|_{L^\infty(\mathbb{R})}, \ \ \ \forall t\leq T.
\end{equation}
\end{proposition}
Proposition \ref{T3.1} can be proved in the standard way, and we omit its proof details. The proof of Proposition \ref{T3.2} is crucial, and is left in next section.

\subsection{The continuity arguments} By the continuity arguments, we can extend the local solution to the global one for all $t\in[0,+\infty)$.  Define 
 \begin{equation}
 \epsilon_*:=\min \left\{\frac{\epsilon_0}{4}, \sqrt{\frac{\epsilon_0^2}{64C_0}-\delta_R^\frac{8}{33}}\right\},\qquad  M:=\frac{\epsilon_0}{4},\notag
\end{equation}
where $\epsilon_0$ and $C_0$ are given in Proposition \ref{T3.2}. By the smallness of $\delta_0$ and $\delta_R\leq\delta_0$ , we have $\frac{\epsilon_0^2}{64C_0}-\delta_R^\frac{8}{33}>0$. Note that $\|\phi_0\|_{H^1(\mathbb{R})}=\|\hat\phi_0\|_{H^1(\mathbb{R})}$. Suppose $\|\hat\phi_0\|_{H^1(\mathbb{R})}< \epsilon_*\leq \frac{\epsilon_0}{4}(=M)$, by local existence result in Proposition \ref{T3.1}, there is a positive constant $T_1=T_1(M)$ such that a unique solution exists on $[0,T_1]$ and satisfies $\|\hat\phi\|_{H^1(\mathbb{R})}\leq 2\|\hat\phi_0\|_{H^1(\mathbb{R})}\leq \frac{\epsilon_0}{2}$ for $t\in[0,T_1]$. 
By Sobolev inequality $\|\hat\phi_0\|_{L^{\infty}(\mathbb{R})} \leq\sqrt 2\| \hat\phi_0\|_{H^1(\mathbb{R})}$, we get $\|u_0\|_{L^{\infty}(\mathbb{R})}<C$, by the maximum principle which implies $\|\phi\|_{L^{\infty}(\mathbb{R})}<C$. Hence, \eqref{3.13}-\eqref{3.15} hold for $t\in[0,T_1]$. Without loss of generality, we can choose $T_1$ small enough such that for all $t\in[0,T_1]$, we have
 \begin{equation*}
  \begin{aligned}
\|u^S &(\cdot-\sigma t+ \mathbf X(t))+u^R(t,\cdot+ \mathbf X(t))-u^S(\cdot-\sigma t)-u^R(t,\cdot)\|_{H^1(\RR)}\\[2mm]
 &\leq Ct\|u^S_\xi+u^R_\xi\|_{H^1(\RR)}\leq \frac{\epsilon_0}{8}.
\end{aligned}
\end{equation*}
Therefore, it holds for $t\in[0,T_1]$ that
 \begin{equation*}
 \begin{aligned}
\di \|\phi(t,\cdot)\|_{H^1(\RR)} \leq &\|u^S(\cdot-\sigma t+ \mathbf X(t))+u^R(t,\cdot+ \mathbf X(t))-u^S(\cdot-\sigma t)-u^R(t,\cdot)\|_{H^1(\RR)}\\
&+ \|\hat\phi(t,\cdot)\|_{H^1(\RR)}\leq \frac{\epsilon_0}{2}+\frac{\epsilon_0}{8}<\epsilon_0.
\end{aligned}
\end{equation*}
Especially, since $ \mathbf X(t)$ is absolutely continuous, and $\hat\phi\in C([0,T_1];H^1{(\mathbb{R})})$, we have $\phi\in C([0,T_1];H^1{(\mathbb{R})})$. We now consider the maximal existence time:
 \begin{equation}
 T_M:=\sup\left\{t>0\Big|\displaystyle \sup_{\tau\in [0,t]}\|\phi(\tau, \cdot)\|_{H^1(\RR)}\leq \epsilon_0\right\}. \notag
\end{equation}
If $T_M<\infty$, then the continuity argument implies that  $\displaystyle \sup_{\tau\in [0,T_M]}\|\phi(\tau, \cdot)\|_{H^1(\mathbb{R})}=\epsilon_0$. It holds from Theorem \ref{T3.2} that
 \begin{equation}
\sup_{\tau\in [0,T_M]}\|\phi(\tau, \cdot)\|_{H^1(\mathbb{R})}\leq \sqrt {C_0(\|\phi_0\|^2_{H^1(\mathbb{R})}+\delta_R^{\frac{8}{33}})}\leq \frac{\epsilon_0}{8},  \notag
\end{equation}
which contradicts $\displaystyle \sup_{\tau\in [0,T_M]}\|\phi(\tau, \cdot)\|_{H^1(\mathbb{R})}=\epsilon_0$.
Thus, $T_M=\infty$, which together with Theorem \ref{T3.2} implies
\begin{equation}\label{3-2}
\begin{array}{ll}
\displaystyle \sup_{t>0}\|\phi(t,\cdot)\|^2_{H^1(\RR)}+\int_0^\infty\|\phi_\xi\|^2_{H^1(\RR)}d\tau+\int_0^\infty\int_\mathbb{R} {\phi}^2((u^S)^{\mb{X}}_\xi+(u^R)^{\mb{X}}_\xi)d \xi d\tau\\[4mm]
\di +\int_0^\infty|  \mathbf{ \dot  X}(t)|^2d\tau\leq C_0(\|\phi_0\|^2_{H^1(\RR)}+\delta_R^\frac{8}{33}),
\end{array}
\end{equation}
and 
\begin{equation}\label{3-3}
| \mathbf{ \dot  X}(t)|\leq C\|\phi(t, \cdot)\|_{L^\infty(\mathbb{R})}, \ \ \ \forall t>0.
\end{equation}

In addition, since the rarefaction wave $u^r$ is Lipschitz continuous in $x$ for all $t>0$ and from Lemma \ref{L2-2}, we have 
\begin{equation}\label{115}
u(t,x)- \Big(u^S\big(x-\sigma t + \mathbf{X}(t);u_-,u_m)+u^r\Big(\frac{x+ \mathbf{X}(t)}{t};u_m,u_+\Big)-u_m\Big)\in C(0,+\infty; H^1(\mathbb R)).
\end{equation}
Since $ \mathbf X(t)$ is absolutely continuous, we have 
\begin{equation}\label{a115}
u(t,x)- \Big(u^S\big(x-\sigma t + \mathbf{X}(t);u_-,u_m)+u^r\Big(\frac{x}{t};u_m,u_+\Big)-u_m\Big)\in C(0,+\infty; H^1(\mathbb R)),
\end{equation}
which is the first line of \eqref{2.21}.

Since $\phi_{\xi\xi}\in  L^2(0,+\infty; L^2(\mathbb R))$ by \eqref{3-2}, and $(u^R)^{\mb{X}}_{\xi\xi} \in  L^2(0,+\infty; L^2(\mathbb R))$ by Lemma \ref{L2-2}, we have 
\begin{equation}
u_{xx}(t,x)- u^S_{xx}(x-\sigma t + \mathbf{X}(t);u_-,u_m) \in L^2(0,+\infty; L^2(\mathbb R)),\notag
\end{equation}
which is the second line of \eqref{2.21}.

In order to get $\|\phi_0\|_{H^1(\RR)}\leq \epsilon_*$, we can choose $\epsilon^*$ and $\delta^*$ which are given in Theorem \ref{T2.1} suitably. Using Lemma \ref{L2-2}, we have 
\begin{align}\label{b2}
&\big{ \|}u^R(1,x;u_m,u_+)-u_m\big{\|}_{L^2(\mathbb{R_-})}\notag\\
&\leq C\delta_R\Big (\int_{-\infty}^{0}e^{-4|x-3u_m^2|}dx\Big)^\frac{1}{2}\leq C\delta_R,
\end{align}
\begin{align}\label{b3}
\di \big{ \|}u^R&(1,x;u_m,u_+)-u_+\big{\|}_{L^2(\mathbb{R_+})} \notag\\
\leq &C\Big (\int_{0}^{3u_+^2}u^R(1,x;u_m,u_+)-u_+)^2dx\Big)^\frac{1}{2} \notag\\
&+C\Big (\int_{3u_+^2}^{+\infty}(u^R(1,x;u_m,u_+)-u_+)^2dx\Big)^\frac{1}{2}\notag\\
\leq& C\Big (\int_{0}^{3u_+^2}(u_m-u_+)^2dx\Big)^\frac{1}{2} +C\delta_R\Big (\int_{3u_+^2}^{+\infty}e^{-4|x-3u_+^2|}dx\Big)^\frac{1}{2} \leq C\delta_R,
\end{align}
and
\begin{equation}\label{b4}
\big{ \|}u^R_x(1,x;u_m,u_+)\big{\|}_{L^2(\mathbb{R})}\leq C\delta_R.
\end{equation}

By \eqref{b2}, \eqref{b3} and \eqref{b4}, we have
\begin{align}\label{b1}
&\|\phi_0\|_{H^1(\RR)}=\big{ \|}u_0-\big(u^S(x;u_-,u_m)+u^R(1,x;u_m,u_+)-u_m\big)\big{\|}_{H^1(\RR)}\notag\\
&=\big{ \|}u_0-\big(u^S(x;u_-,u_m)+u^R(1,x;u_m,u_+)-u_m\big)\big{\|}_{L^2(\mathbb{R_-})}\notag\\
&\quad +\big{ \|}u_0-\big(u^S(x;u_-,u_m)+u^R(1,x;u_m,u_+)-u_m\big)\big{\|}_{L^2(\mathbb{R_+})}\notag\\
&\quad +\big{ \|}u_{0x}(\cdot)-u^S_x(\cdot;u_-,u_m)-u^R_x(1,x;u_m,u_+)\big{\|}_{L^2(\mathbb{R})}\notag\\
&\leq \big{ \|}u_0-u^S(x;u_-,u_m)\big{\|}_{L^2(\mathbb{R_-})}+\big{ \|}u^R(1,x;u_m,u_+)-u_m\big{\|}_{L^2(\mathbb{R_-})}\notag\\
&\quad +\big{ \|}u_0-\big(u^S(x;u_-,u_m)+u_+-u_m\big)\big{\|}_{L^2(\mathbb{R_+})}+\big{ \|}u^R(1,x;u_m,u_+)-u_+\big{\|}_{L^2(\mathbb{R_+})}\notag\\
&\quad +\big{ \|}u_{0x}(\cdot)-u^S_x(\cdot;u_-,u_m)\big{\|}_{L^2(\mathbb{R})}+\big{ \|}u^R_x(1,x;u_m,u_+)\big{\|}_{L^2(\mathbb{R})}\notag\\
&\leq \big{ \|}u_0-u^S(x;u_-,u_m)\big{\|}_{L^2(\mathbb{R_-})}+\big{ \|}u_0-\big(u^S(x;u_-,u_m)+u_+-u_m\big)\big{\|}_{L^2(\mathbb{R_+})}\notag\\
& \di \quad +\big{ \|}u_{0x}(\cdot)-u^S_x(\cdot;u_-,u_m)\big{\|}_{L^2(\mathbb{R})}+C_1\delta_R\notag,
\end{align}
Thus, once we take $\epsilon^*=\frac{\epsilon_*}{2}$ and $\delta^*= \min\left\{\delta_0, \frac{\epsilon_*}{2C_1}\right\}$. Then it holds that $\|\phi_0\|_{H^1(\RR)}\leq \epsilon_*$.\\
 \subsection{Time-asymptotic behaviors} In order to justify the time-asymptotic behaviors \eqref{2.22} and \eqref{2.23}, we set $\di g(t):=\|\phi_\xi(t,\cdot)\|^2$. It is obvious that $g(t)\in L^1(0,+\infty)$. We only need to show that $g^\prime(t)\in L^1(0,+\infty)$. By \eqref{3-2}, \eqref{4.59} and \eqref{4.78}, we have
\begin{align}
\int_0^{+\infty}|g^\prime(t)|dt=&\int_0^{+\infty}\Big|\frac{d}{dt}\|\phi_\xi\|^2\Big|dt\notag\\
\leq &C\int_0^{+\infty}\Big(\|\phi_\xi\|_{H^1(\RR)}^2+|  \mathbf{ \dot  X}(t)|^2+\int_\mathbb{R} {\phi}^2(u^S)^{\mb{X}}_\xi d \xi+\int_\mathbb{R} {\phi}^2(u^R)^{\mb{X}}_\xi d \xi\Big)dt\notag\\
&+C\int_0^{+\infty}\int_\mathbb{R} |F^{\mb{X}}|^2d\xi dt<+\infty.\notag
\end{align}
Hence 
\begin{equation}\label{3.20}
\displaystyle\lim_{t\rightarrow+\infty }g(t)=\displaystyle\lim_{t\rightarrow+\infty} \|\phi_\xi(t,\cdot)\|^2=0.
\end{equation}
By Sobolev inequality, we have 
\begin{equation}\label{3.21}
\displaystyle\lim_{t\rightarrow+\infty}   \|\phi(t, \cdot)\|_{L^{\infty}(\mathbb{R})}\leq\displaystyle\lim_{t\rightarrow+\infty} \sqrt 2\|\phi(t, \cdot)\|^\frac{1}{2}\|\phi_\xi(t, \cdot)\|^\frac{1}{2}=0.         
\end{equation}
In addition, by \eqref{3.19} and  \eqref{3.21}, it holds 
\begin{equation}\label{3.22}
 | \mathbf{ \dot  X}(t)|\leq C\|\phi(t, \cdot)\|_{L^\infty(\mathbb{R})}\rightarrow 0, \ \ \ as \ \ \ t\rightarrow +\infty.
\end{equation}
By \eqref{2.13}, we have $\forall t\geq1$, there exists a uniform constant $C>0$ such that
 \begin{equation}
 \displaystyle \sup_{x\in \mathbb{R}} |w^r(\frac{x}{1+t})-w^r(\frac{x}{t})|\leq C(1+t)^{-1}, \notag
\end{equation}
which results in 
 \begin{equation}\label{sd}
 \displaystyle \sup_{x\in \mathbb{R}} |u^r(\frac{x}{1+t})-u^r(\frac{x}{t})|\leq C(1+t)^{-1}.
\end{equation}
 From Lemma \ref{L2-2} and \eqref{3.22}-\eqref{sd}, it follows that
\begin{equation}\label{3.4}
  \begin{aligned}
 \displaystyle \sup_{x\in \mathbb{R}}| \bar u(t,x)- \tilde u^{\mb{X}}(t,x)|
 &= \displaystyle \sup_{x\in \mathbb{R}}|u^r\big(\frac{x}{t};u_m,u_+\big)-u^R\big(1+t,x+\mathbf{X}(t);u_m,u_+\big)|\\
 &\leq \displaystyle \sup_{x\in \mathbb{R}}\big |u^r\big(\frac{x}{t};u_m,u_+\big)-u^r\big(\frac{x}{1+t};u_m,u_+\big)\big|\\
 &+ \displaystyle \sup_{x\in \mathbb{R}}\Big|u^r\big(\frac{x}{1+t};u_m,u_+\big)- u^r\big(\frac{x+\mathbf{X}(t)}{1+t};u_m,u_+\big)\Big|\\
 &+ \displaystyle \sup_{x\in \mathbb{R}}\Big|u^r\big(\frac{x+\mathbf{X}(t)}{1+t};u_m,u_+\big)-u^R\big(1+t,x+\mathbf{X}(t);u_m,u_+\big)\Big|\\
& \rightarrow 0, (t\rightarrow\infty).
  \end{aligned}
\end{equation}
By \eqref{3.21} and \eqref{3.4}, we prove \eqref{2.22}.   Hence, the proof of Theorem \ref{T2.1} is completed.

\noindent{\bf Notations.} Throughout this paper,  several positive  generic constants are denoted by $C$  which is uniform in time. Define 
  \begin{equation}\label{2.25}
  u_m:=-\frac{u_-}{2}, \ \ \  u_*:=\frac{u_m}{2}=-\frac{u_-}{4}.\notag
   \end{equation}
and  the far field states $u_\pm$ satisfy  \eqref{2.17} means that $u_-<0<u_*< u_m<u_+$.

 $L^p(\mathbb{R})(1\leq p\leq+\infty) $ and $H^1(\mathbb{R})$ denote the usual Lebesgue space and Sobolev space in $\mathbb{R}$ with norm
 \begin{equation}\label{2.26}
\|f\|_{L^p(\mathbb{R})}:=\Big(\int_{\mathbb{R}} |f|^pd\xi\Big)^\frac{1}{p},\quad\|f\|:=\|f\|_{L^2(\mathbb{R})},\ \ \ \|f\|_{H^1(\mathbb{R})}:  =(\|f\|^2+\|f_\xi\|^2)^\frac{1}{2}. \notag
   \end{equation}
   
   For any function $f:\RR^+\times\RR\rightarrow \RR$ and any time-dependent shift $\mathbf{X}(t)$, we denote
    \begin{equation}
f^{\pm\mathbf{X}}(t,\xi):=f(t,\xi\pm\mathbf{X}(t)). \notag
   \end{equation}

 \section{A priori estimates.}  In this section, we prove the uniform-in-time a priori estimates in Proposition \ref{T3.2}. We assume that the Cauchy problem \eqref{3.10} has a solution $ \phi\in C([0,T];H^1(\RR))\cap L^2(0,T; H^2(\RR))$ for some time $T>0$. We start from $L^2$  relative entropy estimate for $\phi$. 
 

 \begin{proposition}\label{P4.1}There exist positive constants $\delta_0,\epsilon_1>0$ such that if the rarefaction wave strength $\delta_R:=|u_+-u_m|\leq\delta_0$ and $ \mathcal{E}(T):=\displaystyle \sup_{0\leq t\leq T}\|\phi(t,\cdot)\|_{H^1(\RR)} \leq\epsilon_1$, then there exists a uniform constant $C>0$ such that for all $t\in[0,T]$
 \begin{align}
\|\phi(t,\cdot)\|^2+&\int_0^t\|\phi_\xi(\tau, \cdot)\|^2d\tau+\int_0^t\int_\mathbb{R} \big(\big|(u^S)^{\mb{X}}_\xi\big|+\big|(u^R)^{\mb{X}}_\xi\big|\big) {\phi}^2d \xi d\tau\\
+& \int_0^t|  \mathbf{ \dot  X}(\tau)|^2d\tau\leq C(\|\phi_0\|^2+\delta_R^\frac{8}{33}).\nonumber
\end{align}
  \end{proposition}

The proof of Proposition \ref{P4.1} can be given by the following Lemmas \ref{L4.2+}  to \ref{L4.11}.
      \begin{lemma}\label{L4.2+} Let $w(u^S(\x))$ be the weight function defined  in \eqref{4.2}.
      Then it holds that
    \begin{equation}\label{4.5}
    \frac{1}{2}\frac{d}{dt}\int_{\mathbb{R}} (\phi^{\mb{-X}})^2 wd\xi+\mathbf{\dot X}(t)\mathbf{Y}(t)+\mathcal{J}^{good}(t)+\mathcal{J}^{bad}(t)=\mathbf{F}(t),
   \end{equation}
   where
    \begin{equation}\label{d4.6}
    \begin{aligned}
\mathbf{Y}(t):=\int_{\mathbb{R}} \phi^{-\mb{X}} w u^S_\xi d\xi+\int_{\mathbb{R}} \phi^{-\mb{X}} w u^R_\xi d\xi-\frac{1}{2}\int_{\mathbb{R}} (\phi^{-\mb{X}})^2 w^\prime u^S_\xi d\xi,
 \end{aligned}
   \end{equation}
   
    \begin{equation}\label{d4.7}
     \begin{aligned}  
   \mathcal{J}^{good}(t)&:=   \int_{\mathbb{R}}\mu w (\phi^{-\mb{X}}_\xi)^2d\xi-3\int_{\mathbb{R}}(\phi^{-\mb{X}})^2 (u^S)^2w^\prime u^S_\xi d\xi\\
   &\quad +3\int_{\mathbb{R}}(\phi^{-\mb{X}})^2(u^R-u_m)w u^S_\xi d\xi +3\int_{\mathbb{R}}(\phi^{-\mb{X}})^2u^Rw u^R_\xi d\xi\\
 & \quad -\frac{3}{2}\int_{\mathbb{R}}(\phi^{-\mb{X}})^2(u^R-u_m)^2w^\prime u^S_\xi d\xi-\frac{3}{4}\int_{\mathbb{R}} (\phi^{-\mb{X}})^4 w^\prime u^S_\xi d\xi,
    \end{aligned}
   \end{equation}
       \begin{equation}\label{d4.8}
     \begin{aligned}  
   \mathcal{J}^{bad}(t):=&3\int_{\mathbb{R}}(\phi^{-\mb{X}})^2wu^Su^S_\xi d\xi+3\int_{\mathbb{R}}(\phi^{-\mb{X}})^2w(u^S-u_m)u^R_\xi d\xi\\
   &-3\int_\mathbb{R}(\phi^{-\mb{X}})^2(u^R-u_m)u^Sw^\prime u^S_\xi d\xi+\sigma\int_{\mathbb{R}}(\phi^{-\mb{X}})^2w^\prime u^S_\xi d\xi \\
   &-\frac{1}{2}\int_{\mathbb{R}}(\phi^{-\mb{X}})^2\mu(u^S_\x)^2w^{\prime\prime}d\x+\int_{\mathbb{R}} (\phi^{-\mb{X}})^3 w (u^S_\xi+u^R_\xi) d\xi\\
   &-2\int_{\mathbb{R}} (\phi^{-\mb{X}})^3  (u^S+u^R-u_m) w^\prime u^S_\xi d\xi,
   \end{aligned}
   \end{equation}
   and
    \begin{equation}\label{d4.9}
\mathbf{F}(t):=-\int_{\mathbb{R}}  F \phi^{-\mb{X}}w d\xi.
   \end{equation}
   \end{lemma}
   
\noindent {\bf Remark 4.1.}	Since $w^\prime\leq 0$, $w^{\prime\prime}\geq 0$, $u^S_\x>0, u^R_\x>0 $ and $ u^R>u_m$,  $\mathcal{J}^{good}(t)$ consists of good terms, while $\mathcal{J}^{bad}(t)$ consists of bad terms and $\mathbf{F}(t) $ can be controlled by wave interaction estimates in  Lemma \ref{L4.9} and the decay properties of approximate rarefaction wave.
\begin{proof}   
  Multiplying \eqref{3.10} by $w^{\mb{X}}\phi$, we can get
 \begin{equation}\label{d4.10}
 \begin{aligned}
 \Big(\frac{1}{2}\phi^2w^{\mb{X}}\Big)_t&-\frac{1}{2}\mathbf{ \dot  X}(t)\phi^2w^{\mb{X}}_\xi-\sigma\phi_\xi\phi w^{\mb{X}}+\Big(f(\phi+\tilde u^{\mb{X}})-f(\tilde u^{\mb{X}} )\Big)_\xi \phi w^{\mb{X}}\\[3mm]
 &+\mathbf{\dot X}(t)\big[(u^S)^{\mb{X}}_\xi+(u^R)^{\mb{X}}_\xi\big]\phi w^{\mb{X}}-\mu \phi_{\xi\xi}\phi w^{\mb{X}}
 =-F^{\mb{X}}\phi w^{\mb{X}}.
 \end{aligned}
\end{equation}
Integrating \eqref{d4.10} over $\RR$ with respect to $\xi$ and changing variable $\xi\rightarrow \xi-\mathbf{  X}(t)$, we have 
 \begin{equation}\label{d4.11}
 \begin{aligned}
    \frac{1}{2}\frac{d}{dt}\int_{\mathbb{R}} (\phi^{-\mb{X}})^2 wd\xi+\mathbf{ \dot  X}(t) \mathbf{Y}(t)+\int_{\mathbb{R}} \Big(f(\phi^{-\mb{X}}+\tilde u)-f(\tilde u)\Big)_\xi \phi^{-\mb{X}} wd\xi\\
    -\sigma\int_{\mathbb{R}} \phi^{-\mb{X}}_\xi\phi^{-\mb{X}} wd\xi- \int_{\mathbb{R}}\mu \phi^{-\mb{X}}_{\xi\xi}\phi^{-\mb{X}} wd\xi=\mathbf{F}(t),
 \end{aligned}
\end{equation}
where $\mathbf{Y}(t)$ and $\mathbf{F}(t)$ are defined in \eqref{d4.6} and \eqref{d4.9} respectively.
Further calculations yield

 \begin{equation}\label{d4.13}
 \begin{aligned}
\int_{\mathbb{R}}-\sigma \phi^{-\mb{X}}_\xi\phi^{-\mb{X}} wd\xi&=-\frac{1}{2}\sigma\int_{\mathbb{R}} (\phi^{-\mb{X}})^2_\xi w d\xi=\frac{\sigma}{2}\int_{\mathbb{R}} (\phi^{-\mb{X}})^2 w^\prime u^S_\xi d\xi,
 \end{aligned}
\end{equation}
and 
 \begin{equation}\label{d4.14}
 \begin{aligned}
 - \int_{\mathbb{R}}\mu \phi^{-\mb{X}}_{\xi\xi}&\phi^{-\mb{X}} wd\xi=\int_{\mathbb{R}}\mu w (\phi^{-\mb{X}}_\xi)^2d\xi-\frac{1}{2}\int_{\mathbb{R}}\mu (\phi^{-\mb{X}})^2w_{\xi\xi}d\xi\\[3mm]
&=\int_{\mathbb{R}}\mu w (\phi^{-\mb{X}}_\xi)^2d\xi+\frac{\sigma}{2}\int_{\mathbb{R}} (\phi^{-\mb{X}})^2 w^\prime u^S_\xi d\xi-\frac{3}{2}\int_\mathbb{R}(\phi^{-\mb{X}})^2w^\prime (u^S)^2u^S_\xi d\xi\\[3mm]
&\di \quad -\frac{1}{2}\int_{\mathbb{R}}\mu (\phi^{-\mb{X}})^2w^{\prime\prime}( u^S_{\xi})^2d\xi.
 \end{aligned}
\end{equation}
It holds that
 \begin{equation}\label{d4.12}
 \begin{aligned}
& \int_{\mathbb{R}} \Big(f(\phi^{-\mb{X}}+\tilde u)-f(\tilde u)\Big)_\xi \phi^{-\mb{X}} wd\xi=-\int_\mathbb{R} (\phi^{-\mb{X}} w)_\xi \Big((\phi^{-\mb{X}}+\tilde u)^3-\tilde u^3\Big)d\xi\\
=&-\int_\mathbb{R} \phi^{-\mb{X}}_\xi\big((\phi^{-\mb{X}})^3+3(\phi^{-\mb{X}})^2\tilde u+3\tilde u^2\phi^{-\mb{X}}\big)w d\xi\\
&-\int_\mathbb{R} w_\xi\big((\phi^{-\mb{X}})^4+3(\phi^{-\mb{X}})^3\tilde u+3(\phi^{-\mb{X}})^2\tilde u^2\big)d\xi \\
=&-\int_\mathbb{R}\Big[\frac{((\phi^{-\mb{X}})^4)_\xi}{4}+((\phi^{-\mb{X}})^3\tilde u)_\xi-(\phi^{-\mb{X}})^3\tilde u_\xi+\big(\frac{3}{2}\tilde u^2(\phi^{-\mb{X}})^2\big)_\xi-3\tilde u\tilde u_\xi(\phi^{-\mb{X}})^2\Big]wd\xi\\
&-\int_\mathbb{R}w_\xi\big[(\phi^{-\mb{X}})^4+3(\phi^{-\mb{X}})^3\tilde u+3(\phi^{-\mb{X}})^2\tilde u^2\big]d\xi \\
=&\int_\mathbb{R} 3\tilde u\tilde u_\xi(\phi^{-\mb{X}})^2w d\xi-\frac{3}{2}\int_\mathbb{R}(\phi^{-\mb{X}})^2\tilde u^2 w_\xi d\xi+\int_\mathbb{R} (\phi^{-\mb{X}})^3\tilde u_\xi wd\xi\\
&-\int _\mathbb{R}\Big(\frac{3}{4}(\phi^{-\mb{X}})^4+2(\phi^{-\mb{X}})^3\tilde u\Big)w_\xi d\xi\\
=&3\int_\mathbb{R} (u^S+u^R-u_m)(u^S_\xi+u^R_\xi)(\phi^{-\mb{X}})^2w d\xi-\frac{3}{2}\int_\mathbb{R}(\phi^{-\mb{X}})^2(u^S+u^R-u_m)^2 w_\xi d\xi\\
&+\int_\mathbb{R} (\phi^{-\mb{X}})^3(u^S_\xi +u^R_\xi )wd\xi-\int _\mathbb{R}\Big(\frac{3}{4}(\phi^{-\mb{X}})^4+2(\phi^{-\mb{X}})^3(u^S+u^R-u_m)\Big) w_\xi d\xi\\
=&3\int_\mathbb{R} (u^S+u^R-u_m)(u^S_\xi+u^R_\xi)(\phi^{-\mb{X}})^2w d\xi-\frac{3}{2}\int_\mathbb{R}(\phi^{-\mb{X}})^2(u^R-u_m)^2w^\prime u^S_\xi d\xi\\
&-\frac{3}{2}\int_\mathbb{R}(\phi^{-\mb{X}})^2(u^S)^2w^\prime u^S_\xi d\xi-3\int_\mathbb{R}(\phi^{-\mb{X}})^2(u^R-u_m)u^Sw^\prime u^S_\xi d\xi\\
&+\int_{\mathbb{R}} (\phi^{-\mb{X}})^3 w (u^S_\xi+u^R_\xi) d\xi-2\int_{\mathbb{R}} (\phi^{-\mb{X}})^3  (u^S+u^R-u_m) w^\prime u^S_\xi d\xi\\
&-\frac{3}{4}\int_{\mathbb{R}} (\phi^{-\mb{X}})^4 w^\prime u^S_\xi d\xi.
 \end{aligned}
\end{equation}

Hence, by \eqref{d4.11}-\eqref{d4.12}, we proved Lemma \ref{L4.2+}.
\end{proof}   

   \begin{lemma}\label{L4.2} It holds that 
    \begin{equation}\label{d4,9}
    \mathbf{\dot X}(t)\mathbf{Y}(t)+\mathcal{J}^{good}(t)+\mathcal{J}^{bad}(t)=\mathbf{G}^S(t)+\mathbf{G}^R(t)+\mathbf{G}^{SR}(t)+\mathbf{N}(t)+\mathbf{J}(t),
   \end{equation}
   where
    \begin{equation}\label{4.6}
    \begin{aligned}
\mathbf{G}^S(t):=&\di \int_{\mathbb{R}}\mu w (\phi^{-\mb{X}}_\xi)^2d\xi+\int_{\mathbb{R}}(\phi^{-\mb{X}})^2 u^S_\xi \Big(\sigma w^\prime-3(u^S)^2w^\prime+3u^Sw-\frac{w^{\prime\prime}}{2}\mu u^S_\xi\Big)d\xi\\
&+\mathbf{ \dot  X}(t)\int_{\mathbb{R}} \phi^{-\mb{X}} w u^S_\xi d\xi-\frac{3}{4}\int_{\mathbb{R}} (\phi^{-\mb{X}})^4 w^\prime u^S_\xi d\xi,
 \end{aligned}
   \end{equation}
       \begin{equation}\label{4.7}
     \begin{aligned}
\mathbf{G}^R(t):=3\int_{\mathbb{R}}(\phi^{-\mb{X}})^2wu^Ru^R_\x d\x,\hspace{7.2cm}
 \end{aligned}
   \end{equation}
     \begin{equation}\label{4.8}
     \begin{aligned}
\mathbf{G}^{SR}(t):=3\int_{\mathbb{R}}(\phi^{-\mb{X}})^2w&(u^R-u_m)u^S_\x d\x-\frac{3}{2}\int_{\mathbb{R}}(\phi^{-\mb{X}})^2(u^R-u_m)^2w^\prime u^S_\xi d\xi,
 \end{aligned}
   \end{equation}
   \begin{equation}\label{4.9}
     \begin{aligned}
\mathbf{N}(t):=&\mathbf{ \dot  X}(t)\int_{\mathbb{R}} \phi^{-\mb{X}} w u^R_\xi d\xi-\frac{1}{2}\mathbf{ \dot  X}(t)\int_{\mathbb{R}} (\phi^{-\mb{X}})^2 w^\prime u^S_\xi d\xi\\
&+\int_{\mathbb{R}} (\phi^{-\mb{X}})^3 w (u^S_\xi+u^R_\xi) d\xi-2\int_{\mathbb{R}} (\phi^{-\mb{X}})^3  (u^S+u^R-u_m) w^\prime u^S_\xi d\xi\hspace{7mm}\\
&-3\int_\mathbb{R}(\phi^{-\mb{X}})^2(u^R-u_m)u^Sw^\prime u^S_\xi d\xi,
 \end{aligned}
   \end{equation}
and
    \begin{equation}\label{4.10}
      \begin{aligned}
\mathbf{J}(t):=3\int_{\mathbb{R}}(\phi^{-\mb{X}})^2(u^S-u_m)wu^R_\x d\x. \hspace{6.5cm}
 \end{aligned}
   \end{equation}
   \end{lemma}
\noindent {\bf Remark 4.2.} First, it is clear that $\mathbf{G}^{R}(t)$ and $\mathbf{G}^{SR}(t)$ are good terms, and $\mathbf{G}^{S}(t)$  are in fact  good terms by Lemma \ref{L4.5} if we choose the weight function $w$ defined in \eqref{4.2} and the time-dependent shift  $\mathbf{X}(t)$ defined in \eqref{3.11}.
 However,   $\mathbf{N}$(t) are bad terms, which can be controlled by Lemma \ref{L4.6} under the assumptions that the perturbations and the rarefaction wave strength are suitably small. While $\mathbf{F}(t)$, $\mathbf{J}(t)$ can be controlled by wave interaction estimates in Lemma \ref{L4.9} and the decay properties of the approximate rarefaction wave.
 
To prove Proposition \ref{P4.1}, we first estimate $\mathbf{G}^S(t),\mathbf{N}(t),\mathbf{J}(t)$ and $\mathbf{F}(t)$ as follows.

\noindent {\bf Estimation of  the $\mathbf{G}^S(t)$}. From \eqref{4.2}, we have $w^\prime=w^{\prime\prime}\equiv 0$ for $\x\in (\xi_*, +\infty)$ and then we can  split $\mathbf{G^S}(t)$ as follows:
\begin{equation}\label{4.6}
\begin{array}{ll}    
\di \mathbf{G^S}(t)=\int_{-\infty}^{\xi_*}\mu w (\phi^{-\mb{X}}_\xi)^2d\xi+\int_{-\infty}^{\xi_*}(\phi^{-\mb{X}})^2 u^S_\xi \Big(\sigma w^\prime-3(u^S)^2w^\prime+3u^Sw-\frac{w^{\prime\prime}}{2}\mu u^S_\xi\Big)d\xi\\[5mm]
\di\qquad +\int_{\xi_*}^{+\infty}\mu w (\phi^{-\mb{X}}_\xi)^2d\xi+3\int_{\xi_*}^{+\infty}(\phi^{-\mb{X}})^2 u^S w u^S_\xi d\xi\\[5mm]
\di \qquad +\mathbf{ \dot  X}(t)\int_{\mathbb{R}} \phi^{-\mb{X}} w u^S_\xi d\xi-\frac{3}{4}\int_{\mathbb{R}} (\phi^{-\mb{X}})^4 w^\prime u^S_\xi d\xi\\[5mm]
\di \quad =\int_{-\infty}^{\xi_*}\mu w (\phi^{-\mb{X}})_\xi^2d\xi+\int_{-\infty}^{\xi_*}(\phi^{-\mb{X}})^2 u^S_\xi \Big(\sigma w^\prime-3(u^S)^2w^\prime+3u^Sw-\frac{w^{\prime\prime}}{2}\mu u^S_\xi\Big)d\xi\\[5mm]
\di\qquad +\frac{8}{25u^2_m}\Big(\int_{-\infty}^{\xi_*} \phi^{-\mb{X}} w u^S_\xi d\xi\Big)^2+\int_{\xi_*}^{+\infty}\mu w (\phi^{-\mb{X}}_\xi)^2d\xi+3\int_{\xi_*}^{+\infty}(\phi^{-\mb{X}})^2 u^S w u^S_\xi d\xi\\[5mm]
\di \qquad+\mathbf{ \dot  X}(t)\int_{\mathbb{R}} \phi^{-\mb{X}} w u^S_\xi d\xi -\frac{8}{25u^2_m}\Big(\int_{-\infty}^{\xi_*} \phi^{-\mb{X}} w u^S_\xi d\xi\Big)^2-\frac{3}{4}\int_{\mathbb{R}} (\phi^{-\mb{X}})^4 w^\prime u^S_\xi d\xi\\[5mm]
\di \quad :=\mathbf {G}^S_1(t)+\mathbf{G}^S_2(t),
\end{array}
\end{equation}
where 
\begin{equation}\label{4.15}
\begin{aligned}
\mathbf{G}^S_1(t):=&\di \int_{-\infty}^{\xi_*}\mu w (\phi^{-\mb{X}}_\xi)^2d\xi+\int_{-\infty}^{\xi_*}(\phi^{-\mb{X}})^2 u^S_\xi \Big(\sigma w^\prime-3(u^S)^2w^\prime+3u^Sw-\frac{w^{\prime\prime}}{2}\mu u^S_\xi\Big)d\xi\\
\qquad&\di  +\frac{8}{25u^2_m}\Big(\int_{-\infty}^{\xi_*} \phi^{-\mb{X}} w u^S_\xi d\xi\Big)^2,
 \end{aligned}
\end{equation}
and
\begin{equation}\label{4.15+}
\begin{aligned}
\mathbf{G}^S_2(t):=&\di \int_{\xi_*}^{+\infty}\mu w (\phi^{-\mb{X}}_\xi)^2d\xi+3\int_{\xi_*}^{+\infty}(\phi^{-\mb{X}})^2 u^S w u^S_\xi d\xi+\mathbf{ \dot  X}(t)\int_{\mathbb{R}} \phi^{-\mb{X}} w u^S_\xi d\xi\\
&\di -\frac{8}{25u^2_m}\Big(\int_{-\infty}^{\xi_*} \phi^{-\mb{X}} w u^S_\xi d\xi\Big)^2-\frac{3}{4}\int_{\mathbb{R}} (\phi^{-\mb{X}})^4 w^\prime u^S_\xi d\xi.
 \end{aligned}
\end{equation}
\begin{lemma}\label{L4.3} $\mathbf{G}^S_1(t)$ defined in \eqref{4.15} satisfies that
\begin{equation}\label{4.16}
\mathbf{G}^S_1(t)\geq \frac{1}{6}\int_{-\infty}^{\xi_*}\mu w (\phi^{-\mb{X}}_\xi)^2d\xi+\frac{4}{5}u_m^3\int_{-\infty}^{\xi_*}(\phi^{-\mb{X}})^2 u^S_\xi d\xi.
\end{equation}
\end{lemma}
\begin{proof}   Let  $y:=\frac{u^S(\x)-u_-}{u_*-u_-}$, then we have $\xi\in(-\infty,\xi_*]\iff y \in[0,1]$. Since $y_\x=\frac{2}{5u_m}u^S_\x>0$,  there exists a unique inverse function $\x=\x(y)$ by the inverse function theorem. For any fixed $t>0$, we can write  $\psi(t,y):= \phi^{-\mb{X}}(t,\x) w(u^S(\x))$.  
In order to use the weighted Poincar${\acute{\rm e}}$ inequality with $\psi(t,y)$, we first have
\begin{equation}\label{4.17}
 \begin{aligned}
2\Big(\int_0^1\psi(t,y) dy\Big)^2&=2\Big(\int_{-\infty }^{\xi_*} \phi^{-\mb{X}} w \frac{u^S_\xi}{u_*-u_-}d\xi\Big)^2=\frac{8}{25u^2_m}\Big(\int_{-\infty}^{\xi_*} \phi^{-\mb{X}} w u^S_\xi d\xi\Big)^2.
 \end{aligned}
\end{equation}
Furthermore, we have 
\begin{equation}\label{4.18}
2\int_0^1\psi^2(t,y)dy=2\int_{-\infty }^{\xi_*} (\phi^{-\mb{X}})^2w^2\frac{u^S_\xi}{u_*-u_-}d\xi=\frac{4}{5u_m}\int_{-\infty }^{\xi_*}( \phi^{-\mb{X}})^2w^2u^S_\xi d\x.
\end{equation}
and 
\begin{equation}\label{4.19}
 \begin{aligned}
&\di \int_0^1 y(1-y)  |\psi_y(t,y)|^2dy 
=\int_{-\infty }^{\xi_*} {[(\phi^{-\mb{X}} w)_\x]}^2\frac{(u^S-u_-)(u_*-u^S)}{u^S_\xi}d\xi\frac{1}{u_*-u_-} \\[4mm]
&\quad \di =\frac{2}{5u_m}\int_{-\infty }^{\xi_*} {[(\phi^{-\mb{X}} w)_\x]}^2\frac{\mu(u_*-u^S)}{(u_m-u^S)^2}d\xi\\[4mm]
&\quad\di =\frac{2}{5u_m}\int_{-\infty }^{\xi_*} (\phi^{-\mb{X}}_\xi)^2 w^2\frac{\mu(u_*-u^S)}{(u_m-u^S)^2}d\xi
+\frac{2}{5u_m}\int_{-\infty }^{\xi_*} (\phi^{-\mb{X}})^2 w_\xi^2\frac{\mu(u_*-u^S)}{(u_m-u^S)^2}d\xi\\[4mm]
&\di\qquad +\frac{4}{5u_m}\int_{-\infty }^{\xi_*} \phi^{-\mb{X}} \phi^{-\mb{X}}_\xi ww_\xi\frac{\mu(u_*-u^S)}{(u_m-u^S)^2}d\xi:=I_1+I_2+I_3.
 \end{aligned}
\end{equation}
Noticing that 
\begin{equation}\label{4.20}
 \begin{aligned}
 \mu w_{\x\x}&=\mu[w^{\prime}(u^S)u^S_\x]_\x=\mu w^{\prime\prime}(u^S)(u^S_\x)^2+\mu w^\prime(u^S)u^S_{\x\x}
 \\
 &=(u^S-u_-)(u^S-u_m)^2w^{\prime\prime}(u^S)u^S_\x+[-3u_m^2+3(u^S)^2]w^\prime(u^S)u^S_\x,
 \end{aligned}\notag
\end{equation} we have
\begin{equation}\label{4.20}
 \begin{aligned}
I_3&\di =\frac{4}{5u_m}\int_{-\infty }^{\xi_*} \phi^{-\mb{X}} \phi^{-\mb{X}}_\xi ww_\xi\frac{\mu(u_*-u^S)}{(u_m-u^S)^2}d\xi  = \frac{2}{5u_m}\int_{-\infty }^{\xi_*} ((\phi^{-\mb{X}})^2)_\xi ww_\xi\frac{\mu(u_*-u^S)}{(u_m-u^S)^2}d\xi \\[2mm]
&\di=-\frac{2}{5u_m}\int_{-\infty }^{\xi_*} (\phi^{-\mb{X}})^2  \Big(ww_\xi\frac{\mu(u_*-u^S)}{(u_m-u^S)^2}\Big)_\xi d\xi \\[2mm]
&\di =-I_2-\frac{2}{5u_m}\int_{-\infty }^{\xi_*} (\phi^{-\mb{X}})^2ww_{\xi\xi}\frac{\mu(u_*-u^S)}{(u_m-u^S)^2}d\xi \\[2mm]
 &\di\quad-\frac{2}{5u_m}\int_{-\infty }^{\xi_*} (\phi^{-\mb{X}})^2 ww_\xi \Big(\frac{\mu(u_*-u^S)}{(u_m-u^S)^2}\Big)_\xi d\xi \\[2mm]
&\di 
=-I_2+\frac{2}{5u_m}\int_{-\infty }^{\xi_*} (\phi^{-\mb{X}})^2ww^\prime u^S_\xi(2u^S+2u_m-u_*)d\xi\\[2mm]
&\di \quad -\frac{2}{5u_m}\int_{-\infty }^{\xi_*} (\phi^{-\mb{X}})^2ww^{\prime\prime}(u_*-u^S)(u^S+2u_m)u^S_\xi d\xi.
\end{aligned}
\end{equation}
By \eqref{4.19} and \eqref{4.20}, we have 
\begin{equation}\label{4.21}
 \begin{aligned}
\di  \int_0^1y(1-y)  |\psi_y(t,y)|^2dy&=I_1+I_2+I_3\\
&\di =\frac{2}{5u_m}\int_{-\infty }^{\xi_*} (\phi^{-\mb{X}}_\xi)^2 w^2\frac{\mu(u_*-u^S)}{(u_m-u^S)^2}d\xi\\[3mm]
&\di \quad+\frac{2}{5u_m}\int_{-\infty }^{\xi_*} (\phi^{-\mb{X}})^2ww^\prime u^S_\xi(2u^S+2u_m-u_*)d\xi\\[3mm]
&\di \quad
-\frac{2}{5u_m}\int_{-\infty }^{\xi_*}( \phi^{-\mb{X}})^2ww^{\prime\prime}(u_*-u^S)(u^S+2u_m)u^S_\xi d\xi.
 \end{aligned}
\end{equation}
By weighted Poincar${\acute{\rm e}}$ inequality \eqref{4.1+}, and \eqref{4.17}, \eqref{4.18}, \eqref{4.21}, it holds that
\begin{equation}\label{4.22}
 \begin{aligned}
&\di  \mathbf{G}^S_1(t)\geq  \int_{-\infty}^{\xi_*}\mu w (\phi^{-\mb{X}}_\xi)^2 \Big(1-w\frac{(u_*-u^S)}{(u_m-u^S)^2}\frac{2}{5u_m}\Big)d\xi\\[2mm]
&\di \  +\frac{2}{5u_m}\int_{-\infty}^{\xi_*}(\phi^{-\mb{X}})^2 u^S_\xi\underbrace{ \Big(\sigma w^\prime-3(u^S)^2w^\prime+3u^Sw-\frac{w^{\prime\prime}}{2}\mu u^S_\xi\Big)(u_*-u_-)}_{H_1}d\xi\\[2mm]
&\di\   +\frac{2}{5u_m}\int_{-\infty}^{\xi_*}(\phi^{-\mb{X}})^2 u^S_\xi \underbrace{w\Big( w^{\prime\prime}(u_*-u^S)(u^S+2u_m)+2w-w^\prime (2u^S+2u_m-u_*)\Big)}_{H_2}d\xi.
 \end{aligned}
\end{equation}
By \eqref{4.2}, for $ \xi\in(-\infty,\xi_*),$ we  have
\begin{equation}\label{4.23}
1-w\frac{(u_*-u^S)}{(u_m-u^S)^2}\frac{2}{5u_m}=
\begin {cases}
\di 1-\frac{\frac{u_m}{2}-u^S}{u_m-u^S}=\frac{\frac{u_m}{2}}{u_m-u^S}, &\di  \xi\in(-\infty,\xi_1),\\[5mm]
\di \frac{1}{2u_m^3(u_m-u^S)}[8(u^S)^4-4u_m(u^S)^3+u_m^4], &\di \xi\in [\xi_1,\xi_*).\\
\end{cases}
\end{equation}
For $\x\in (-\infty, \xi_1)$, or equivalently, $u^S\in(u_-,0),$ we have
\begin{equation}\label{4.23+}
\frac{\frac{u_m}{2}}{u_m-u^S}> \frac{\frac{u_m}{2}}{u_m-u_-}=\frac{\frac{u_m}{2}}{u_m-(-2u_m)}=\frac 16.
\end{equation}
For $\x\in [\xi_1,\xi_*)$, or equivalently, $u^S\in[0,\frac{u_m}{2}),$ we have
\begin{equation}\label{4.24}
 \begin{aligned}
 \frac{1}{2u_m^3(u_m-u^S)}[8(u^S)^4-4u_m(u^S)^3+u_m^4]>\frac{1}{4},
 \end{aligned}
\end{equation}
 since if we set $g(u^S)=8(u^S)^4-4u_m(u^S)^3+u_m^4$, then $g^\prime(u^S)=(u^S)^2(32u^S-12u_m)$, and it is easy to check that  for
$u^S\in [0, \frac{u_m}{2})$, $g_{\min}(u^S)=g(\frac{3u_m}{8})=\frac{485}{512}u_m^4\approx0.947u_m^4$  and $\frac{1}{2u_m^3(u_m-u^S)}\geq\frac{1}{2u_m^4}$. 
 
By \eqref{4.23}-\eqref{4.24}, we have for $ \xi\in(-\infty,\xi_*),$
\begin{equation}\label{231}
1-w\frac{(u_*-u^S)}{(u_m-u^S)^2}\frac{2}{5u_m}> \frac 16.
\end{equation}
By the definition of the weight function $w$ defined in \eqref{4.2},  we can compute 
\begin{equation}\label{4.25}
H_1+H_2=
\begin {cases}
\di \frac{25}{8}u_m^3(u_m-u^S),\qquad \qquad  \qquad \qquad \qquad \qquad  \qquad \di \mbox{for}\ \xi\in(-\infty,\xi_1),\\[4mm]
\di -\frac{1800}{u_m^4}(u^S)^8+\frac{400}{u_m^3}(u^S)^7+\frac{3950}{u_m^2}(u^S)^6-\frac{3225}{u_m}(u^S)^5+\frac{775}{2}(u^S)^4\\[5mm]
\di \quad+325u_m(u^S)^3-\frac{75}{2}u_m^2(u^S)^2-\frac{25}{8}u_m^3u^S+\frac{25}{8}u_m^4, \ \ \di \mbox{for}\ \xi \in [\xi_1,\xi_*),\\
\end{cases}
\end{equation}
where $H_1$ and $H_2$ are defined in \eqref{4.22}. The detailed calculations can be found in Section 5.

Now we claim that for both $\xi \in (-\infty, \x_1)$ and $\xi\in[\xi_1,\xi_*)$,
\begin{equation}\label{claim}
H_1+H_2 > 2u_m^4.
\end{equation}
First, it is obvious that for $\xi \in (-\infty, \x_1)$, $u^S\in (u_-, 0)$,   and then $\di \frac{25}{8}u_m^3(u_m-u^S)>\frac{25}{8}u_m^4>2u_m^4$. For $\xi\in[\xi_1,\xi_*)$, or equivalently, $u^S\in[0,\frac{u_m}{2})$, define
\begin{equation}
 \begin{aligned}
 h(u^S):=&\di h_1(u^S)-\frac{3225}{u_m}(u^S)^5+\frac{775}{2}(u^S)^4+325u_m(u^S)^3\\[3mm]
&\di -\frac{75}{2}u_m^2(u^S)^2-\frac{25}{8}u_m^3u^S+\frac{25}{8}u_m^4,
\end{aligned}\notag
\end{equation}
where 
$$
\begin{array}{ll}
h_1(u^S)& \di :=-\frac{1800}{u_m^4}(u^S)^8+\frac{400}{u_m^3}(u^S)^7+\frac{3950}{u_m^2}(u^S)^6\\[5mm]
&\di =\frac{(u^S)^6}{u_m^2}\left[-\frac{1800}{u_m^2}(u^S)^2+\frac{400}{u_m}u^S+3950\right]\geq  \frac{3700}{u_m^2}(u^S)^6
\end{array}
$$
for $u^S\in[0,\frac{u_m}{2})$.
Then we  have
\begin{equation}
 \begin{aligned}
h(u^S)\geq&\underbrace{\frac{3700}{u_m^2}(u^S)^6-\frac{3225}{u_m}(u^S)^5+703(u^S)^4}_{h_2(u^S)}\\
&-\frac{631}{2}(u^S)^4+325u_m(u^S)^3-\frac{75}{2}u_m^2(u^S)^2-\frac{25}{8}u_m^3u^S+\frac{25}{8}u_m^4\\[4mm]
\geq&-\frac{631}{2}(u^S)^4+325u_m(u^S)^3-\frac{75}{2}u_m^2(u^S)^2-\frac{25}{8}u_m^3u^S+\frac{25}{8}u_m^4\\[4mm]
\geq&\underbrace{\frac{669}{4}u_m(u^S)^3-\frac{75}{2}u_m^2(u^S)^2-\frac{25}{8}u_m^3u^S+\frac{25}{8}u_m^4}_{l(u^S)}, 
\end{aligned}\notag
\end{equation}
where we have used the fact that $h_2(u^S)\geq 0$, $\forall u^S\in[0,\frac{u_m}{2})$. We can compute that $$l^\prime(u^S)=\frac{2007}{4}u_m(u^S)^2-75u_m^2u^S-\frac{25}{8}u_m^3,$$
which have two roots $\di u_1=\frac{150-2\sqrt{11896.875}}{2007}u_m<0$ and $\di u_2=\frac{150+2\sqrt{11896.875}}{2007}u_m\in (0, u_*)$. Thus for $u^S\in[0,u_2),l^\prime(u^S)<0$ and for $u^S\in(u_2,u_*),l^\prime(u^S)>0.$ Therefore, we have $l_{\min}(u^S)=l(u_1)>2u_m^4$, and then the claim in \eqref{claim} is proved. 
Thus the proof of Lemma \ref{L4.3} is completed.
 \end{proof} 
\begin{lemma}\label{L4.4} It holds that
\begin{equation}\label{4.26}
\begin{array}{ll}
\di \mathbf{G}^S_2(t)\geq &\di  \int_{\xi_*}^{+\infty}\mu w (\phi_\xi^{-\mb{X}})^2d\xi+(\frac{45}{16}-\frac{9}{8}\ln 2)u_m^3\int_{\xi_*}^{+\infty}(\phi^{-\mb{X}})^2 u^S_\xi d\xi\\[5mm]
&\di +\frac{25}{64}u_m^2|\mathbf{ \dot  X}(t)|^2+\frac{3}{4}\int_{\mathbb{R}} (\phi^{-\mb{X}})^4 |w^\prime| u^S_\xi d\xi,
\end{array}
\end{equation}
where $\mathbf{G}^S_2(t)$ is defined by \eqref{4.15} and $\mathbf{ \dot  X}(t)$ is defined in \eqref{3.11}.
\end{lemma}
\begin{proof} It is sufficient to show that 
\begin{equation}
 \begin{aligned}
&\di 3\int_{\xi_*}^{+\infty}(\phi^{-\mb{X}})^2 u^S w u^S_\xi d\xi+\mathbf{ \dot  X}(t)\int_{\mathbb{R}} \phi^{-\mb{X}} w u^S_\xi d\xi-\frac{8}{25u^2_m}\Big(\int_{-\infty}^{\xi_*} \phi^{-\mb{X}} w u^S_\xi d\xi\Big)^2\\
&\di \geq(\frac{45}{16}-\frac{9}{8}\ln 2)u_m^3\int_{\xi_*}^{+\infty}(\phi^{-\mb{X}})^2 u^S_\xi d\xi+\frac{25}{64}u_m^2|\mathbf{ \dot  X}(t)|^2. \notag
\end {aligned}
\end{equation}
Noticing that $\di \mathbf{ \dot  X}(t) =\frac{32}{25u_m^2} \int_{\mathbb{R}} \phi w^{\mb{X}}(u^S)^{\mb{X}}_\xi d\xi=\frac{32}{25u_m^2} \int_{\mathbb{R}} \phi^{-\mb{X}} wu^S_\xi d\xi$, we have
\begin{equation}\label{4.27}
 \begin{aligned}
\frac{1}{2}\mathbf{ \dot  X}(t)  \int_{\mathbb{R}} \phi^{-\mb{X}} wu^S_\xi d\xi&\di =\frac{16}{25u_m^2} \Big( \int_{\mathbb{R}} \phi^{-\mb{X}} wu^S_\xi d\xi\Big)^2\\[2mm]
&\di =\frac{16}{25u_m^2}\Big(\int_{-\infty }^{\xi_*}   \phi^{-\mb{X}} wu^S_\xi d\xi\Big)^2+\frac{16}{25u_m^2} \Big( \int_{\xi_*} ^{+\infty} \phi^{-\mb{X}} wu^S_\xi d\xi\Big)^2\\[2mm]
&\di \qquad +\frac{32}{25u_m^2}\Big(\int_{-\infty }^{\xi_*}   \phi^{-\mb{X}} wu^S_\xi d\xi\Big) \Big(\int_{\xi_*} ^{+\infty}  \phi^{-\mb{X}} wu^S_\xi d\xi\Big)\\[2mm]
&\di \geq \frac{8}{25u_m^2}\Big(\int_{-\infty }^{\xi_*}   \phi^{-\mb{X}} wu^S_\xi d\xi\Big)^2-\frac{16}{25u_m^2} \Big( \int_{\xi_*} ^{+\infty} \phi^{-\mb{X}} wu^S_\xi d\xi\Big)^2,
 \end{aligned}
\end{equation}
where in the last inequality we have used the fact $\di |ab|\leq \frac{a^2}{4}+b^2$ such that
$$
\begin{array}{ll}
\di \frac{32}{25u_m^2}\Big(\int_{-\infty }^{\xi_*}   \phi^{-\mb{X}} wu^S_\xi d\xi\Big) \Big(\int_{\xi_*} ^{+\infty}  \phi^{-\mb{X}} wu^S_\xi d\xi\Big)\\[5mm]
\di  \geq -\frac{8}{25u_m^2}\Big(\int_{-\infty }^{\xi_*}   \phi^{-\mb{X}} wu^S_\xi d\xi\Big)^2-\frac{32}{25u_m^2} \Big( \int_{\xi_*} ^{+\infty} \phi^{-\mb{X}} wu^S_\xi d\xi\Big)^2.
\end{array}
$$
H${ \ddot{\rm o}}$lder inequality yields
$$
\begin{array}{ll}
\di  \Big( \int_{\xi_*} ^{+\infty} \phi^{-\mb{X}} wu^S_\xi d\xi\Big)^2&\di \leq \Big(\int_{\xi_*} ^{+\infty}  3(\phi^{-\mb{X}})^2 u^S  wu^S_\xi d\xi \Big)\Big(\int_{\xi_*} ^{+\infty}   \frac{ wu^S_\xi}{3 u^S} d\xi \Big)\\[5mm]
&\di = \Big(\int_{\xi_*} ^{+\infty}  3(\phi^{-\mb{X}})^2 u^S  wu^S_\xi d\xi \Big)\frac{5 u_m^2}{8}\int_{\xi_*}^{+\infty} \big(\ln|u^S(\xi)|\big)_\x d\xi\\[5mm]
&\di = \Big(\int_{\xi_*} ^{+\infty}  3(\phi^{-\mb{X}})^2 u^S  wu^S_\xi d\xi \Big)\frac{5 u_m^2}{8}\big(\ln|u^S(\xi)|\big)\big|_{\xi_*} ^{+\infty}\\[5mm]
&\di =\frac{5 u_m^2\ln 2}{8}\Big(\int_{\xi_*} ^{+\infty}  3(\phi^{-\mb{X}})^2 u^S  wu^S_\xi d\xi \Big).
\end{array}
$$
Substituting the above inequality into \eqref{4.27} gives that
\begin{equation}\label{4.27+}
 \begin{aligned}
\frac{1}{2}\mathbf{ \dot  X}(t)  \int_{\mathbb{R}} \phi^{-\mb{X}} wu^s_\xi d\xi
&\di \geq\frac{8}{25u_m^2}\Big(\int_{-\infty }^{\xi_*}   \phi^{-\mb{X}} wu^S_\xi d\xi\Big)^2-\frac{2}{5}\ln 2\Big(\int_{\xi_*} ^{+\infty}  3(\phi^{-\mb{X}})^2 u^S  wu^S_\xi d\xi \Big).
 \end{aligned}
\end{equation}
From \eqref{4.27+} and  the fact $u^S(\x)\geq\frac{u_m}{2}, \forall \x \geq \x_*$, one has
\begin{equation}\label{4.28}
 \begin{aligned}
&\di  3\int_{\xi_*}^{+\infty}(\phi^{-\mb{X}})^2 u^S w u^S_\xi d\xi+\mathbf{ \dot  X}(t)\int_{\mathbb{R}} \phi^{-\mb{X}} w u^S_\xi d\xi-\frac{8}{25u^2_m}\Big(\int_{-\infty}^{\xi_*} \phi^{-\mb{X}} w u^S_\xi d\xi\Big)^2\\[3mm]
&\di \geq3(1-\frac{2}{5}\ln 2)\int_{\xi_*}^{+\infty}(\phi^{-\mb{X}})^2 u^S w u^S_\xi d\xi+\frac{1}{2}\mathbf{ \dot  X}(t)\int_{\mathbb{R}} \phi^{-\mb{X}} w u^S_\xi d\xi\\[3mm]
&\di \geq \frac{45}{8}(1-\frac{2}{5}\ln 2)u_m^2\int_{\xi_*}^{+\infty}(\phi^{-\mb{X}})^2 u^S u^S_\xi d\xi+\frac{1}{2}\mathbf{ \dot  X}(t)\int_{\mathbb{R}} \phi^{-\mb{X}} w u^S_\xi d\xi\\[3mm]
  &\di \geq(\frac{45}{16}-\frac{9}{8}\ln 2)u_m^3\int_{\xi_*}^{+\infty}(\phi^{-\mb{X}})^2 u^S_\xi d\xi+\frac{25}{64}u_m^2|\mathbf{ \dot  X}(t)|^2.
 \end{aligned}
\end{equation}
where in the last equality we have used the fact $\di \frac{1}{2}\mathbf{ \dot  X}(t)  \int_{\mathbb{R}} \phi^{-\mb{X}} wu^s_\xi d\xi=\frac{25}{64}u_m^2|\mathbf{ \dot  X}(t)|^2$. 
Thus, the proof of Lemma \ref{L4.4} is complete.
 \end{proof} 
 By  Lemma \ref{L4.3} and  Lemma \ref{L4.4}, we can get the  following lemma.
 \begin{lemma}\label{L4.5}  The following estimation of $\mathbf{G^S}(t)$ holds true.
 \begin{equation}\label{4.29}
 \begin{array}{ll}
 \mathbf{G^S}(t)\geq &\di \frac{5}{16}u_m^2\int_{\mathbb{R}} \mu (\phi_\xi^{-\mb{X}})^2d\xi+\frac{4}{5}u_m^3\int_{\mathbb{R}}  (\phi^{-\mb{X}})^2 u^S_\xi d\xi\\[5mm]
 &\di+\frac{25}{64}u_m^2|\mathbf{ \dot  X}(t)|^2+\frac{3}{4}\int_{\mathbb{R}} (\phi^{-\mb{X}})^4 |w^\prime| u^S_\xi d\xi.
 \end{array}
 \end{equation}
  \end{lemma}

\noindent {\bf Estimation of  the $\mathbf{N}(t)$}. The estimation of  the $\mathbf{N}(t)$ is given in the following lemma.
 \begin{lemma}\label{L4.6} There exists a positive constant $\epsilon_2>0$, such that if $\mathcal{E}(T):=\displaystyle \sup_{0\leq t\leq T}\|\phi\|_{H^1(\mathbb{R})}\leq \epsilon_2$, then it holds that
  \begin{equation}\label{4.30}
  \begin{aligned}
| \mathbf{N}(t)|\leq&  (36\sqrt2+\frac{35\sqrt2}{2})u_m^2\epsilon_2\int_{\mathbb{R}}( \phi^{-\mb{X}})^2 u^S_\xi d\xi\\
 &+\big(5\sqrt 2u_m \delta_R\epsilon_2+15u_m^2\delta_R\big)\int_\mathbb{R} (\phi^{-\mb{X}})^2 u^S_\xi d\xi\\
& +(\sqrt 2\epsilon_2+\frac{48}{5}\delta_R)\int_\mathbb{R}( \phi^{-\mb{X}})^2  u^R_\xi wd\xi+\frac{25}{128}u_m^2|\mathbf{ \dot  X}(t)|^2.
 \end{aligned}
 \end{equation}
 \end{lemma}
 \begin{proof} We can  rewrite  $\mathbf{N}(t)$ as:
  \begin{equation}\label{4.30d}
  \begin{aligned}
 \mathbf{N}(t)=&\underbrace{\int_\mathbb{R}  \big((\phi^{-\mb{X}})^3u^S_\xi w-2(\phi^{-\mb{X}})^3u^Su^S_\xi w^\prime \big) d\xi}_{ \mathbf{N}_1(t)}
 +\underbrace{\dot{\mb{X}}(t)\Big(-\frac{1}{2}\int_\mathbb{R} w^\prime u^S_\xi (\phi^{-\mb{X}})^2d\xi\Big)}_{ \mathbf{N}_2(t)}\\
 &+\underbrace{ \int_\mathbb{R} (\phi^{-\mb{X}})^3 u^R_\xi wd\xi+\mathbf{ \dot  X}(t)\int_{\mathbb{R}} \phi^{-\mb{X}} w u^R_\xi d\xi}_{ \mathbf{N}_3(t)}\\
&\underbrace{-2\int_\mathbb{R} (\phi^{-\mb{X}})^3 (u^R-u_m) w_\xi d\xi-3\int_\mathbb{R} (\phi^{-\mb{X}})^2 (u^R-u_m)u^S w_\xi d\xi}_{ \mathbf{N}_4(t)}.\\
 \end{aligned}
 \end{equation}
 Since $w^\prime\leq 0$ and $u^S_\x>0$, we have 
\begin{equation}\label{4.31}
 \begin{aligned}
  |\mathbf{N}_1(t)|&=\Big|\int_\mathbb{R}  \big((\phi^{-\mb{X}})^3u^S_\xi w-2(\phi^{-\mb{X}})^3u^Su^S_\xi w^\prime\big) d\xi\Big| \\
&=\Big |\int_\mathbb{R} ( \phi^{-\mb{X}})^3\big( w-2 u^S w^\prime\big)u^S_\xi d\xi\Big|\ \\
&\leq\|\phi^{-\mb{X}}\|_{L^\infty (\mathbb{R})}\Big[\|w\|_{L^\infty(\mathbb{R})}+2\|u^S\|_{L^\infty(\mathbb{R})}\|w^\prime\|_{L^\infty(\mathbb{R})} \Big] \int_\mathbb{R}  (\phi^{-\mb{X}})^2 u^S_\xi d\xi\\
&=\Big(\|\phi^{-\mb{X}}\|_{L^\infty(\mathbb{R})}\frac{35}{2}u_m^2\Big) \int_\mathbb{R}  (\phi^{-\mb{X}})^2 u^S_\xi d\xi\\
&\leq\Big(\frac{35\sqrt2}{2}u_m^2\|\phi^{-\mb{X}}\|_{H^1(\mathbb{R})}\Big)\int_\mathbb{R} ( \phi^{-\mb{X}})^2 u^S_\xi d\xi,
 \end{aligned}
\end{equation}
and
\begin{equation}\label{4.32}
 \begin{aligned}
 &| \mathbf{N}_2(t)|=\Big|\dot{\mb{X}}(t)\Big(-\frac{1}{2}\int_\mathbb{R} w^\prime u^S_\xi (\phi^{-\mb{X}})^2d\xi\Big)\Big|\\
&=\Big|\frac{32}{25u_m^2}\Big(\int_{\mathbb{R}} \phi^{-\mb{X}} w u^S_\xi d\xi\Big) \Big(-\frac{1}{2}\int_\mathbb{R}w^\prime u^S_\xi (\phi^{-\mb{X}})^2d\xi\Big)\Big|\\
&\leq\frac{16}{25u_m^2}\|\phi^{-\mb{X}}\|_{L^\infty(\mathbb{R})}\|w\|_{L^\infty(\mathbb{R})}\|w^\prime\|_{L^\infty(\mathbb{R})}\int_{\RR} u^S_\x d\x\ \int_\mathbb{R} (\phi^{-\mb{X}})^2 u^S_\xi d\xi\\
&\leq \frac{48}{25u_m}\|\phi^{-\mb{X}}\|_{L^\infty(\mathbb{R})}\|w\|_{L^\infty(\mathbb{R})}\|w^\prime\|_{L^\infty(\mathbb{R})}\int_\mathbb{R} (\phi^{-\mb{X}})^2 u^S_\xi d\xi\\
&\leq 36\sqrt2 u_m^2 \|\phi^{-\mb{X}}\|_{H^1(\mathbb{R})}\int_\mathbb{R}( \phi^{-\mb{X}})^2 u^S_\xi d\xi.
 \end{aligned}
\end{equation}
  By $\epsilon$-Cauchy inequality, we have $\forall \epsilon>0,$
  $$
  \left|\mathbf{ \dot  X}(t)\int_{\mathbb{R}} \phi^{-\mb{X}} w u^R_\xi d\xi\right|\leq \epsilon |\mathbf{ \dot  X}(t)|^2+\frac{1}{4\epsilon}\Big(\int_{\mathbb{R}} \phi^{-\mb{X}} w u^R_\xi d\xi\Big)^2.
  $$
  Taking $\di \epsilon=\frac{25}{128}u_m^2$ in the above inequality and using H${ \ddot{\rm o}}$lder  inequality yield
      \begin{equation}\label{4.35}
 \begin{aligned}
|\mathbf{N}_3(t)|&=\Big|\int_\mathbb{R} (\phi^{-\mb{X}})^3w u^R_\xi d\xi+\mathbf{ \dot  X}(t)\int_{\mathbb{R}} \phi^{-\mb{X}} w u^R_\xi d\xi\Big|\\
&\leq \|\phi^{-\mb{X}}\|_{L^\infty(\mathbb{R})}\int_\mathbb{R} (\phi^{-\mb{X}})^2  wu^R_\xi d\xi\\
&\di \quad +\frac{25}{128}u_m^2|\mathbf{ \dot  X}(t)|^2+\frac{32}{25u_m^2}\Big( \int_\mathbb{R} \phi^{-\mb{X}} w u^R_\xi d\xi\Big)^2\\
&\leq \|\phi^{-\mb{X}}\|_{L^\infty(\mathbb{R})}\int_\mathbb{R} (\phi^{-\mb{X}})^2 w u^R_\xi d\xi+\frac{25}{128}u_m^2|\mathbf{ \dot  X}(t)|^2\\
&\di \quad +\frac{32}{25u_m^2}\int_\mathbb{R} (\phi^{-\mb{X}})^2 w u^R_\xi d\xi\int_\mathbb{R}  wu^R_\xi d\xi\\
&\leq\Big(\sqrt 2\|\phi^{-\mb{X}}\|_{H^1(\mathbb{R})}+\frac{48}{5}\delta_R\Big)\int_\mathbb{R}( \phi^{-\mb{X}})^2  u^R_\xi wd\xi+\frac{25}{128}u_m^2|\mathbf{ \dot  X}(t)|^2,
   \end{aligned}
    \end{equation}
    where in the last inequality we have used the fact $\di \int_\mathbb{R}  u^R_\xi wd\xi\leq \frac{15}{2}u_m^2  \int_\mathbb{R}  u^R_\xi d\x=\frac{15}{2}u_m^2\delta_R.$
Finally, since $-\frac{5}{2}u_m\leq w^\prime\leq 0,$ we have
  \begin{equation}\label{4.36}
 \begin{aligned}
|\mathbf{N}_4(t)|&=\Big|2\int_\mathbb{R}( \phi^{-\mb{X}})^3 (u^R-u_m) w_\xi d\xi+3\int_\mathbb{R} (\phi^{-\mb{X}})^2 (u^R-u_m)u^S w_\xi d\xi\Big|\\
&\leq \Big[2\delta_R\|\phi^{-\mb{X}}\|_{L^\infty(\mathbb{R})}\|w^\prime\|_{L^\infty(\mathbb{R})}+3\delta_R|u_-|\|w^\prime\|_{L^\infty(\mathbb{R})}\Big]\int_\mathbb{R} (\phi^{-\mb{X}})^2 u^S_\xi d\xi\\
&\leq\big(5\sqrt 2u_m \delta_R\|\phi^{-\mb{X}}\|_{H^1(\mathbb{R})}+15u_m^2\delta_R\big)\int_\mathbb{R} (\phi^{-\mb{X}})^2 u^S_\xi d\xi.
 \end{aligned}
    \end{equation}
By \eqref{4.30d}-\eqref{4.36}, we proved Lemma \ref{L4.6}.
 \end{proof} 
Integrating \eqref{4.5} with respect to $t$ and using Lemmas \ref{L4.2}-\ref{L4.6}, we have the following lemma.
 \begin{lemma}\label{L4.8}There exist positive constants $\delta_0>0$ and $\epsilon_3>0$, such that if $\mathcal{E}(T):=\displaystyle \sup_{0\leq t\leq T}\|\phi\|_{H^1(\RR)}\leq\epsilon_3$ and the rarefaction wave strength $\delta_R\leq\delta_0$, then  it holds that $\forall t\in [0, T],$
   \begin{equation}\label{4.37}
    \begin{aligned}
   \frac{1}{2}\int_\mathbb{R}(&\phi^{-\mb{X}})^2wd\xi+ \frac{5}{16}u_m^2\int_0^t \int_\mathbb{R}\mu(\phi_\xi^{-\mb{X}})^2d\xi d\tau+\frac{u_m^3}{2}\int_0^t \int_\mathbb{R}(\phi^{-\mb{X}})^2 (u^S_\xi+u^R_\xi) d\xi d\tau \\
 &  +\frac{25}{128}u_m^2\int_0^t |\mathbf{ \dot  X}(\tau)|^2 d\tau\leq  \frac{1}{2}\int_\mathbb{R}\phi_0^2wd\xi+\int_0^t  \mathbf{F}(\tau)d\tau-\int_0^t  \mathbf{J}(\tau)d\tau.
     \end{aligned}
    \end{equation}
 \end{lemma}
It remains to estimate the terms $\di \int_0^t  \mathbf{F}(\tau) d\tau$ and $\di \int_0^t  \mathbf{J}(\tau)d\tau$. First, we handle the wave interaction estimates. Note that the shifted Oleinik shock and the shifted rarefaction wave are always attached together, therefore, their wave interactions are very subtle.

   \begin{lemma}\label{L4.9} (Wave interaction estimates)  There exists generic uniform-in-time constant $C>0$, such that
 \begin{equation}\label{4.38}
 \int_{-\infty}^{0} |u^S-u_m|u^R_\xi d\xi\leq C \delta_S \delta_R^{\frac{2}{11}} (1+t)^{-\frac{4}{5}} ,
 \end{equation}
  \begin{equation}\label{4.39}
   \int^{+\infty}_{0} |u^S-u_m|u^R_\xi d\xi\leq C\delta_S^{-\frac{3}{5}}\delta_R^{\frac{1}{5}}(1+t)^{-\frac{4}{5}},
   \end{equation}
     \begin{equation}\label{4.40}
 \int_{-\infty}^{0} |u^R-u_m|u^S_\xi d\xi\leq C \delta_S \delta_R^{\frac{2}{11}} (1+t)^{-\frac{4}{5}} ,
 \end{equation}
     \begin{equation}\label{4.41}
\int_{0}^{[f^\prime(u_+)-\sigma](t+1)}|u^R-u^r|u^S_\xi d\xi \leq C \delta_R^{\frac{1}{5}}\delta_S (1+t)^{-\frac{4}{5}},
 \end{equation}
      \begin{equation}\label{4.42}
\int_{0}^{[f^\prime(u_+)-\sigma](t+1)}|u^r-u_m|u^S_\xi d\xi \leq C \delta_R^{\frac{1}{5}}\delta_S^{-\frac{3}{5}}(1+t)^{-\frac{4}{5}} \ln^{\frac{4}{5}}(1+C\delta_S^2\delta_R t),
 \end{equation}
   \begin{equation}\label{4.43}
   \int_{[f^\prime(u_+)-\sigma](t+1)}^{+\infty}|u^R-u_m|u^S_\xi d\xi \leq C  \delta_R \delta_S(1+C \delta_R \delta_S^2t)^{-1},
    \end{equation}
where  $u^R=u^R\big(1+t,\xi+\sigma t;u_m,u_+\big)$ and $u^r=u^r\big(\frac{\xi+\sigma (1+t)}{1+t};u_m,u_+\big)$.
\end{lemma}
 \begin{proof}
 By Lemma \ref{L2-1} and Lemma \ref{L2-2}, $\forall \ve\in (0,1),$
 \begin{equation}\label{4.44}
 \begin{array}{ll}
\di \int_{-\infty}^{0} |u^S-u_m|u^R_\xi d\xi&\di \leq\delta_S \int_{-\infty}^{0} u^R_\xi d\xi\\
&\di =\delta_S|u^R(1+t,\sigma t)-u_m|\leq C_\ve \delta_S \delta_R^{\frac{2\ve}{2+\ve}}(1+t)^{-1+\ve},
 \end{array}
\end{equation}
which is exactly \eqref{4.38} if we take $\ve=\frac 15$.
By H${\ddot{\rm o}}$lder inequality, $\forall\ve\in(0,1)$ 
     \begin{equation}\label{4.45}
 \begin{aligned}
 \int^{+\infty}_{0} |u^S-u_m|u^R_\xi d\xi\leq& \Big( \int^{+\infty}_{0} |u^S-u_m|^{\frac{1}{1-\ve}}d\xi \Big)^{1-\ve}\Big( \int^{+\infty}_{0}|u^R_\xi|^\frac{1}{\ve}d\xi \Big)^{\ve}\\
& \leq  C\delta_S  \Big(\int^{+\infty}_{0}(1+C\delta_S^2\xi)^{-\frac{1}{1-\ve}}d\xi\Big)^{1-\ve}\delta_R^{\ve}(1+t)^{-1+\ve}\\
& =C_\ve \delta_S^{2\ve-1}\delta_R^{\ve}(1+t)^{-1+\ve} \Big(\int^{+\infty}_{0}(1+\xi)^{-\frac{1}{1-\ve}}d\xi\Big)^{1-\ve}\\
&\leq C\delta_S^{-\frac{3}{5}}\delta_R^{\frac{1}{5}}(1+t)^{-\frac{4}{5}}.
 \end{aligned}
    \end{equation}
If we take $\ve=\frac 15$ in\eqref{4.45}, then \eqref{4.39} is proved.
By Lemma \ref{L2-2} (5),  we have for $\forall \ve\in (0,1),$
   \begin{equation}\label{4.46}
 \begin{aligned}
    \int_{-\infty}^{0} |u^R-u_m|u^S_\xi d\xi&\di \leq C_\ve(1+t)^{-1+\ve}\delta_R^{\frac{2\ve}{2+\ve}}  \int_{-\infty}^{0} u^S_\xi d\xi\\
    &\di \leq C_\ve(1+t)^{-1+\ve}\delta_R^{\frac{2\ve}{2+\ve}} \delta_S,
 \end{aligned}
    \end{equation}
which is exactly \eqref{4.40} if we take $\ve=\frac 15$.  
Now we prove \eqref{4.41}.  It holds that
       \begin{equation}\label{4.47}
 \begin{aligned}
    &\di  \int_{0}^{[f^\prime(u_+)-\sigma](1+t)}|u^R-u^r|u^S_\xi d\xi\\
      & \leq \int_{0}^{[f^\prime(u_+)-\sigma](1+t)}|u^R(1+t,\xi+\sigma t)-u^R(1+t,\xi+\sigma (1+t))|u^S_\xi d\xi \\
   &  \quad+\int_{0}^{[f^\prime(u_+)-\sigma](1+t)}|u^R(1+t,\xi+\sigma (1+t))-u^r|u^S_\xi d\xi\\[1mm]
    &\leq  C\|u^R_{\xi}\|_{L^{\infty}(\mathbb{R})}\int_{0}^{+\infty}u^S_\xi d\xi + C_\ve(1+t)^{-1+\ve}\delta_R^{\ve}   \int_{0}^{+\infty}u^S_\xi d\xi  \\
    &\leq C_\ve(1+t)^{-1+\ve}\delta_R^{\ve} \delta_S,
 \end{aligned}
    \end{equation}
    which is  \eqref{4.41} if we take $\ve=\frac 15$.  
    
  On the other hand, noticing that $ C^{-1}\delta_R\leq f^\prime(u_+)-\sigma\leq C\delta_R$ for some constant $C>0$ and the inequality $\ln(1+x)\leq x$ for $x>0$, we have
    \begin{equation}\label{4.48}
 \begin{aligned}
&\di  \int_{0}^{[f^\prime(u_+)-\sigma](1+t)}|u^r-u_m|u^S_\xi d\xi\\
 &= \int_{0}^{[f^\prime(u_+)-\sigma](1+t)}\Big|(f^\prime)^{-1}(\frac {\xi+\sigma (1+t)}{1+t}\Big)-(f^\prime)^{-1}(\sigma)\Big|u^S_\xi d\xi\\
    &\leq C\delta_S ^3  \int_{0}^{[f^\prime(u_+)-\sigma](1+t)} \frac{\xi}{1+t}\frac{1}{(C\delta_S^2\xi+1)^2}d\xi \\
  &\leq C{\delta_S } \frac{1}{1+t} \int_{0}^{[f^\prime(u_+)-\sigma](1+t)} \frac{1}{C\delta_S^2\xi+1}d\xi\\
  &\leq C\delta_S ^{-1} \frac{\ln[C\delta_S^2\delta_R (1+t)+1]}{1+t} \leq  C\delta_S ^{-\frac{3}{5}}\delta_R^{\frac{1}{5}}\Big( \frac{\ln(C\delta_S^2\delta_R t+1)}{1+t}\Big)^\frac{4}{5},
        \end{aligned}
    \end{equation}
    and 
    \begin{equation}\label{4.49}
 \begin{aligned}
& \int_{[f^\prime(u_+)-\sigma](1+t)}^{+\infty}|u^R-u_m|u^S_\xi d\xi\\
 &\di  \leq \delta_R\Big|u_m-u^S\big([f^\prime(u_+)-\sigma](1+t)\big)\Big| \\
  &\leq C\delta_R\delta_S\frac{1}{1+C\delta_S^2\big([f^\prime(u_+)-\sigma]t\big)}\leq C  \delta_R \delta_S\big(1+C \delta_R \delta_S^2t\big)^{-1}.
   \end{aligned}
    \end{equation}
    Thus, the proof of Lemma \ref{L4.9} is complete.
\end{proof}

In the following estimation, $C$  may depend on $\delta_S$ but is independent of  $\delta_R$. Since $\delta_R$ is small and we can assume $\delta_R\leq1$. Thus we have the following lemmas.

\noindent{\bf Estimation of  the $\di \int_0^t   \mathbf{J}(\tau)d\tau$}.  
 \begin{lemma}\label{L4.10}  $ \mathbf{J}(t)$ is defined in \eqref{4.10}, it holds that
  \begin{equation}\label{4.50}
\Big | \int_0^t \mathbf{J}(\tau) d \tau\Big| \leq \frac{u_m^2}{16} \int_0^t \mu \| \phi^{-\mb{X}}_\xi \|^2 d\tau+C u_m^2  \sup_{0\leq t \leq T}\| \phi^{-\mb{X}} \|^2 \delta_R^{\frac{4}{11}},\ for\ all \ t\in[0,T].
  \end{equation}
 \end{lemma}
 \begin{proof}By Cauchy inequality and Sobolev inequality, using Lemma \ref{L4.9}, we have $\forall t\in [0, T],$
  \begin{equation}\label{4.51}
 \begin{aligned}
 \Big| \int_0^t   \mathbf{J}(\tau) d \tau\Big|&=\Big|3\int_0^t   \int_\mathbb{R} (\phi^{-\mb{X}})^2(u^S-u_m)u^R_\xi wd\xi d \tau\Big|\\
& \leq C u_m^2 \int_0^t \|\phi^{-\mb{X}}\|^2_{L^{\infty}(\mathbb{R})}  \int_\mathbb{R} |u^S-u_m|u^R_\xi d\xi d \tau\\
&\leq  \frac{u_m^2}{16} \int_0^t \mu \| \phi^{-\mb{X}}_\xi \|^2 d\tau+C u_m^2  \sup_{0\leq t \leq T}\| \phi^{-\mb{X}} \|^2  \int_0^t  \Big(\int_\mathbb{R} |u^S-u_m|u^R_\xi d\xi\Big)^2 d \tau\\
& \leq \frac{u_m^2}{16} \int_0^t \mu \| \phi^{-\mb{X}}_\xi \|^2 d\tau+C u_m^2  \sup_{0\leq t \leq T}\| \phi ^{-\mb{X}}\|^2 \delta_R^{\frac{4}{11}},\ for\ all\ t\in[0,T].
   \end{aligned}
    \end{equation}
 \end{proof}
\noindent {\bf Estimation of  the $\di \int_0^t  \mathbf{F}(\tau) d\tau$}. 
  \begin{lemma}\label{L4.11} $\mathbf{F}(t)$ is defined in \eqref{4.9}, it holds that $\forall t\in[0,T]$
  \begin{equation}\label{4.52}
\Big |\int_0^t  \mathbf{F}(\tau) d\tau\Big| \leq \frac{u_m^2}{16} \int_0^t \mu \| \phi^{-\mb{X}}_\xi \|^2 d\tau+C u_m^2  \sup_{0\leq t \leq T}\| \phi^{-\mb{X}} \|^\frac{2}{3} \delta_R^{\frac{8}{33}}.
  \end{equation}
 \end{lemma}
  \begin{proof}

By Cauchy inequality and Sobolev inequality, we have
  \begin{equation}\label{4.53}
 \begin{aligned}
| \mathbf{F}(\tau)|=\Big| \int_\mathbb{R} F\phi^{-\mb{X}} wd\xi\Big|&
\leq Cu_m^2\|\phi^{-\mb{X}}\|_{L^{\infty}(\mathbb{R})} \int_\mathbb{R}|F|d\xi 
\leq Cu_m^2\|\phi^{-\mb{X}}\|^{\frac{1}{2}}\|\phi^{-\mb{X}}_\xi\|^{\frac{1}{2}}\int_\mathbb{R} |F|d\xi\\
&\leq \frac{u_m^2}{16}\mu  \|\phi^{-\mb{X}}_\xi\|^2+Cu_m^2 \|\phi^{-\mb{X}}\|^{\frac{2}{3}}\Big(\int_\mathbb{R} |F|d\xi\Big)^{\frac{4}{3}},
   \end{aligned}
    \end{equation}
where 
\begin{equation}\label{4.54}
  \int_\mathbb{R} |F|d\xi \leq C\int_\mathbb{R} |u^S-u_m|u^R_\xi d\xi+ C\int_\mathbb{R} |u^R-u_m|u^S_\xi d\xi+\mu \int_\mathbb{R} |u^R_{\xi\xi}|d\xi .
    \end{equation}
Using Lemma \ref{L4.9}, we have
\begin{equation}\label{4.55}
 \int_0^{t} \Big(\int_\mathbb{R} |u^S-u_m|u^R_\xi d\xi\Big)^\frac{4}{3}d\tau \leq C\delta_R^{\frac{8}{33}}.
  \end{equation}
  Based on the fact
    \begin{equation}\label{4-56}
 \begin{aligned}
    \int_\mathbb{R} |u^R-u_m|u^S_\xi d\xi\leq  \int_{-\infty}^{0} |u^R-u_m|u^S_\xi d\xi
    + \int_{0}^{[f^\prime(u_+)-\sigma](1+t)}|u^R-u^r|u^S_\xi d\xi\\
    +\int_{0}^{[f^\prime(u_+)-\sigma](1+t)}|u^r-u_m|u^S_\xi d\xi +  \int_{[f^\prime(u_+)-\sigma](1+t)}^{+\infty}|u^R-u_m|u^S_\xi d\xi ,
    \end {aligned}
      \end{equation}
and 
\begin{equation}\label{h}
 \begin{aligned}
 \int_0^{t} \Big( \int_{[f^\prime(u_+)-\sigma](1+t)}^{+\infty}|u^R-u_m|u^S_\xi d\xi\Big)^\frac{4}{3}d\tau \leq C\delta_R^{\frac{4}{3}}\int_0^{+\infty}(1+C\delta_Rt)^{-\frac{4}{3}}dt\leq C\delta_R^{\frac{1}{3}}.
    \end {aligned}\notag
      \end{equation}
Thus, by  Lemma \ref{L4.9}, we obtain
  \begin{equation}\label{4.56}
 \int_0^{t} \Big(\int_\mathbb{R} |u^R-u_m|u^S_\xi d\xi\Big)^\frac{4}{3}d\tau \leq C\delta_R^{\frac{8}{33}}.
  \end{equation}
  By Lemma \ref{2.2}, we have 
  \begin{equation}\label{4.57}
    \int_\mathbb{R} |u^R_{\xi\xi}|d\xi\leq
\left\{\begin{array}{ll}
\delta_R, &\di  1+t\leq \delta_R^{-1},\\[2mm]
(1+t)^{-1}, &\di 1+t\geq \delta_R^{-1},
  \end{array}\right.
    \end{equation}
  thus
  \begin{equation}\label{4.58}
 \int_0^{t} \Big(\int_\mathbb{R} |u^R_{\xi\xi}|d\xi \Big)^\frac{4}{3}d\tau \leq C\delta_R^\frac{1}{3}.
  \end{equation}
Thus we complete the proof of Lemma \ref{L4.11}.
 \end{proof}
 
Notice that $\|\phi\|_{H^1(\mathbb{R})}=\|\phi^{-\mb{X}}\|_{H^1(\mathbb{R})}$ and $\displaystyle \sup_{0\leq t\leq T}\|\phi\|_{H^1(\mathbb{R})}=\displaystyle \sup_{0\leq t\leq T}\|\phi^{-\mb{X}}\|_{H^1(\mathbb{R})}$. Using Lemma \ref{L4.8}, \ref{L4.10} and  \ref{L4.11} and changing of variable $\xi\rightarrow \xi+\mathbf{X}(t)$, we can prove  Proposition \ref{P4.1}.
 
To accomplish the a priori estimates in Proposition \ref{T3.2}, we need to obtain the following $L^2$-estimate for $\phi_\xi$.
  \begin{proposition}\label{P4.2} 
  There exist positive constants $\delta_0,\epsilon_4>0$, if $ \mathcal{E}(T):=\displaystyle \sup_{0\leq t\leq T}\|\phi\|_{H^1(\mathbb{R})}\leq\epsilon_4$ and the rarefaction wave strength $\delta_R=|u_+-u_m|\leq\delta_0$, then there exists a constant $C$ such that for all $t\in[0,T]$, it holds that
 \begin{equation}
\|\phi_\xi (t)\|^2+\int_0^t\|\phi_{\xi\xi}\|^2d\tau\leq C(\|\phi_0\|_{H^1(\mathbb{R})}^2+\delta_R^\frac{8}{33}).\notag
\end{equation}
    \end {proposition}
  \begin{proof} Multiplying \eqref{3.10} by $-\phi_{\xi\xi}$ and integrating the resulted equation with respect to $\xi$, we can obtain
   \begin{equation}\label{4.59}
   \begin{array}{ll}
 \di  \frac{1}{2}\frac{d}{dt}\| \phi_\xi\|^2+\mu\|\phi_{\xi\xi}\|^2 =\int_{\mathbb{R}}\Big(f(\phi+\tilde u^{\mb{X}})-f(\tilde u^{\mb{X}})\Big)_\xi \phi_{\xi\xi}d\xi
 \\[5mm]
 \di\quad\quad+\mathbf{ \dot X}(t)\int_{\mathbb{R}}\phi_{\xi\xi}((u^S)^{\mb{X}}_\xi+(u^R)^{\mb{X}}_\xi)d\xi+\int_{\mathbb{R}}\phi_{\xi\xi} F^{\mb{X}}d\xi.
 \end{array}
  \end{equation}
By Cauchy inequality on the right-hand side of \eqref{4.59}, we obtain
   \begin{equation}\label{4.60}
  \frac{1}{2}\frac{d}{dt}\| \phi_\xi\|^2+\frac{\mu}{2}\|\phi_{\xi\xi}\|^2\leq C\int_{\mathbb{R}}\Big|\Big(f(\phi+\tilde u^{\mb{X}})-f(\tilde u^{\mb{X}})\Big)_\xi\Big|^2d\xi+C\int_{\mathbb{R}}|F^{\mb{X}}|^2d\xi+C|\mathbf{ \dot X}(t)|^2.
  \end{equation}
  The first term on the right-hand side of \eqref{4.60} can be estimated as 
    \begin{equation}\label{4.61}
 \begin{aligned}
 \int_{\mathbb{R}}&\Big|\Big(f(\phi+\tilde u^{\mb{X}})-f(\tilde u^{\mb{X}})\Big)_\xi\Big|^2d\xi \\
& \leq \int_{\mathbb{R}} \Big|\Big(f^\prime(\phi+\tilde u^{\mb{X}})-f^\prime(\tilde u^{\mb{X}})\Big) \tilde u^{\mb{X}}_\xi\Big|^2d\xi+ \int_{\mathbb{R}} \Big|f^\prime(\phi+\tilde u^{\mb{X}})\phi_\xi \Big|^2 d\xi\\
 &\leq C \int_{\mathbb{R}}\phi^2 ((u^S)^{\mb{X}}_\xi)^2d\xi+ C \int_{\mathbb{R}}\phi^2 ((u^R)^{\mb{X}}_\xi)^2d\xi + C\|\phi_\xi\|^2\\
& \leq C\int_{\mathbb{R}}\phi^2 (u^S)^{\mb{X}}_\xi d\xi+C\int_{\mathbb{R}}\phi^2 (u^R)^{\mb{X}}_\xi d\xi+ C\|\phi_\xi\|^2.
   \end{aligned}
    \end{equation}
By Proposition \ref{P4.1}, we have
  \begin{equation}\label{4.62}
  \int_0^t \int_{\mathbb{R}}\Big|\Big(f(\phi+\tilde u^{\mb{X}})-f(\tilde u^{\mb{X}})\Big)_\xi\Big|^2d\xi d\tau\leq C (\|\phi_0\|^2+\delta_R^{\frac{8}{33}}), \ \ \ t\in[0,T].
     \end{equation}
We next estimate the second term on the right-hand side of  \eqref{4.60}. By the definition of $F^{\mb{X}}$, \eqref{3.7}, it holds that

  \begin{equation}\label{4.63}
  \begin{aligned}\int_{\mathbb{R}}|F^{\mb{X}}|^2d\xi=\int_{\mathbb{R}}|F|^2d\xi\leq &C  \int_{\mathbb{R}}|u^S-u_m|^2(u^R_\xi)^2d\xi\\
&+C  \int_{\mathbb{R}}|u^R-u_m|^2(u^S_\xi)^2d\xi+C\int_{\mathbb{R}} |u^R_{\xi\xi}|^2 d\xi.
\end{aligned}
     \end{equation}
By \eqref{4.38}, \eqref{4.39} and Lemma \ref{L2-2}, we have
    \begin{equation}\label{4.64}
 \begin{aligned}
 \int_{\mathbb{R}}&|u^S-u_m|^2(u^R_\xi)^2d\xi\\[3mm]
& \leq C\|u^R_\xi\|_{L^{\infty}(\mathbb R)}\int_{\mathbb{R}}|u^S-u_m|u^R_\xi d\xi\\[3mm]
& \leq C\delta_{R}^\frac{3}{5}(1+t)^{-\frac{2}{5}}\delta_{R}^\frac{2}{11}(1+t)^{-\frac{4}{5}}=C\delta_{R}^\frac{43}{55}(1+t)^{-\frac{6}{5}},
   \end{aligned}
    \end{equation}
which implies that
 \begin{equation}\label{4.65}
   \int_0^t  \int_{\mathbb{R}}|u^S-u_m|^2(u^R_\xi)^2d\xi d\tau\leq C\delta_{R}^\frac{43}{55}.
\end{equation} 
We next estimate the second term on the right-hand side of  \eqref{4.63}.
  \begin{equation}\label{4.66}
 \begin{aligned}
 & \int_{\mathbb{R}}|u^R-u_m|^2(u^S_\xi)^2d\xi\\[4mm]
  & \leq  \int_{-\infty}^0 |u^R-u_m|^2(u^S_\xi)^2d\xi +\int_{0}^{[f^\prime(u_+)-\sigma](1+t)} |u^R-u_m|^2(u^S_\xi)^2d\xi \\[4mm]
& \quad +\int^{+\infty}_{[f^\prime(u_+)-\sigma](1+t)} |u^R-u_m|^2(u^S_\xi)^2d\xi .
    \end{aligned}
    \end{equation}
    For $\di \int_{-\infty}^0 |u^R-u_m|^2(u^S_\xi)^2d\xi $,  by Lemma \ref{L2-2}, we have 
      \begin{equation}\label{4.67}
 \begin{aligned}
   \int_{-\infty}^0 |u^R-u_m|^2(u^S_\xi)^2d\xi   \leq C_{\epsilon} \delta_R^{\frac{4\ve}{2+\epsilon}}(1+t)^{-2+2\epsilon}.
    \end{aligned}
    \end{equation}
Taking $\epsilon=\frac{1}{3}$, we have
    \begin{equation}\label{4.68}
     \int_0^t   \int_{-\infty}^0 |u^R-u_m|^2(u^S_\xi)^2d\xi d\tau\leq C  \delta_R^{\frac{4}{7}}.
      \end{equation}
 For $\di \int_{0}^{[f^\prime(u_+)-\sigma](1+t)} |u^R-u_m|^2(u^S_\xi)^2d\xi$, we have
   \begin{equation}\label{4.69}
 \begin{aligned}
 &  \int_{0}^{[f^\prime(u_+)-\sigma](1+t)}|u^R-u_m|^2(u^S_\xi)^2d\xi\\[2mm]
  & \leq \int_{0}^{[f^\prime(u_+)-\sigma](1+t)} |u^R-u^r|^2(u^S_\xi)^2d\xi+\int_{0}^{[f^\prime(u_+)-\sigma](1+t)} |u^r-u_m|^2(u^S_\xi)^2d\xi.
    \end{aligned}
    \end{equation}
By Lemma \ref{L2-2}, we have 
\begin{equation}\label{4.70}
 \begin{aligned}
& \int_{0}^{[f^\prime(u_+)-\sigma](1+t)} |u^R-u^r|^2(u^S_\xi)^2d\xi\\[2mm]
&\leq\int_{0}^{[f^\prime(u_+)-\sigma](1+t)} |u^R(1+t,\x+\sigma t)-u^R(1+t,\x+\sigma (1+t))|^2(u^S_\xi)^2d\xi\\
&+ \int_{0}^{[f^\prime(u_+)-\sigma](1+t)} |u^R(1+t,\x+\sigma (1+t))-u^r|^2(u^S_\xi)^2d\xi\\[2mm]
&\leq C\|u^R_\x \|_{L^\infty(\RR)}^2 \int_{0}^{+\infty}(u^S_\xi)^2d\xi+C_{\ve} \delta_R^{2\ve}(1+t)^{-2+2\ve} \int_{0}^{+\infty}(u^S_\xi)^2d\xi\\
& \leq C_{\ve} \delta_R^{2\ve}(1+t)^{-2+2\ve},
    \end{aligned}
    \end{equation}
and
\begin{equation}\label{4.71}
 \begin{aligned}
&\int_{0}^{[f^\prime(u_+)-\sigma](1+t)} |u^r-u_m|^2(u^S_\xi)^2d\xi\\[2mm]
&=\int_{0}^{[f^\prime(u_+)-\sigma](1+t)}\Big|\Big(f^\prime)^{-1}(\frac {\xi+\sigma (1+t)}{1+t}\Big)-(f^\prime)^{-1}(\sigma)\Big|^2(u^S_\xi)^2 d\xi\\[2mm]
    &\leq C \int_{0}^{[f^\prime(u_+)-\sigma](1+t)} \frac{\xi^2}{(1+t)^2}\frac{{\delta_S}^6}{(C{\delta_S}^2\xi+1)^4}d\xi \\[2mm]
  &\leq C \frac{1}{(1+t)^2} \int_{0}^{[f^\prime(u_+)-\sigma](1+t)} \frac{1}{C\xi+1}d\xi \\[2mm]
  &\leq C \frac{1}{(1+t)^2} \ln(1+C\delta_Rt)\\[2mm]
  &\leq C\delta_R^{\frac{2}{3}}\frac{\ln^{\frac{1}{3}}(1+C\delta_R t)}{(1+t)^{\frac{4}{3}}}.
    \end{aligned}
    \end{equation}
Taking $\ve=\frac{1}{3}$ in \eqref {4.70} and using \eqref{4.70}, \eqref{4.71}, we have 
\begin{equation}\label{4.72}
 \int_{0}^t  \int_{0}^{[f^\prime(u_+)-\sigma](1+t)} |u^R-u_m|^2(u^S_\xi)^2d\xi d\tau \leq C \delta_R^{\frac{2}{3}}.
 \end{equation}
  For $\di  \int^{+\infty}_{[f^\prime(u_+)-\sigma](1+t)} |u^R-u_m|^2(u^S_\xi)^2d\xi$, we have
  \begin{equation}\label{4.73}
 \begin{aligned}
&\int^{+\infty}_{[f^\prime(u_+)-\sigma](1+t)} |u^R-u_m|^2(u^S_\xi)^2d\xi\\[2mm]
&\leq \delta_R^2\|u^S_\xi\|_{L^{\infty}\big([f^\prime(u_+)-\sigma](1+t),+\infty \big)}\int^{+\infty}_{[f^\prime(u_+)-\sigma](1+t)} u^S_\xi d\xi\\[2mm]
&\leq C\delta_R^2\frac{1}{(1+C\delta_Rt)^2},
    \end{aligned}
    \end{equation}
  which implies 
    \begin{equation}\label{4.74}
   \int_0^t \int^{+\infty}_{[f^\prime(u_+)-\sigma](1+t)} |u^R-u_m|^2(u^S_\xi)^2d\xi d\tau\leq C \delta_R.
      \end{equation}
Substituting \eqref{4.68}, \eqref{4.72}, \eqref{4.74} into  \eqref{4.66}, we have
  \begin{equation}\label{4.75}
   \int_0^t  \int_{\mathbb{R}}|u^R-u_m|^2(u^S_\xi)^2d\xi d\tau\leq C \delta_R ^{\frac{4}{7}}.
      \end{equation}
      We finally estimate the third  term on the right-hand side of  \eqref{4.60}. 
      \begin{equation}\label{4.76}
     \int_{\mathbb{R}}|u^R_{\xi\xi}|^2d\xi \leq C \|u^R_{\xi\xi}\|_{L^{\infty}(\mathbb{R})}^{\frac{2}{3}} \|u^R_{\xi\xi}\|_{L^{\infty}(\mathbb{R})}^{\frac{1}{3}}(1+t)^{-1}\leq C \delta_R ^{\frac{2}{3}}(1+t)^{-\frac43}, 
        \end{equation}
     which implies that
      \begin{equation}\label{4.77}
       \int_0^t \int_{\mathbb{R}}|u^R_{\xi\xi}|^2d\xi d\tau \leq C \delta_R ^{\frac{2}{3}}.
        \end{equation}
        Substituting \eqref{4.65}, \eqref{4.75}, \eqref{4.77} into  \eqref{4.63}, we have
     \begin{equation}\label{4.78}
       \int_0^t \int_{\mathbb{R}}|F^{\mb{X}}|^2d\xi d\tau \leq C \delta_R ^{\frac{4}{7}}.
        \end{equation}   
        Integrating \eqref{4.59} with respect to $t$, by \eqref{4.62},\eqref{4.78} and Proposition \ref{P4.1}, we can prove   Proposition \ref{P4.2}.
        
By Proposition \ref{P4.1} and \ref{P4.2}, we finally have the desired uniform-in-time a priori estimates Proposition \ref{T3.2}.
  \end{proof}
 \section{Appendix}
 
In this appendix we show the detailed calculations of $H_1+H_2$ in \eqref{4.25}. Notice that  $H_1$ and  $H_2$ are defined in \eqref{4.22} with
$$ 
\begin{array}{ll}
\di H_1&\di =\big[(3u_m^2-3(u^S)^2)w^\prime+3u^Sw-\frac{1}{2}w^{\prime\prime}\mu u^S_\xi\big](u_*-u_-), \\ [4mm]
\di H_2&\di =ww^{\prime\prime}(u_*-u^S)(u^S+2u_m)+2w^2+ww^\prime(u_*+u_--2u^S). 
\end{array} 
$$
 When $\xi\in(-\infty,\xi_1)$, we have $w=\frac{5}{2}u_m(u_m-u^S),\ w^\prime=-\frac{5}{2}u_m,\ w^{\prime\prime}=0$. Then
  \begin{equation}\label{a18}
 \begin{aligned}
H_1&=\frac{5}{2}u_m\big[(3u_m^2-3(u^S)^2)\cdot(-\frac{5}{2}u_m)+3u^S\cdot\frac{5}{2}u_m(u_m-u^S)\big]\\[3mm]
&=\frac{25}{4}u_m^2\big[3(u^S)^2-3u_m^2+3u_mu^S-3(u^S)^2\big]\\[3mm]
&=-\frac{75}{4}u_m^4+\frac{75}{4}u_m^3u^S,
  \end{aligned}
   \end{equation}  
   and 
  \begin{equation}\label{a19}
 \begin{aligned}  
H_2&=2[\frac{5}{2}u_m(u_m-u^S)]^2+\frac{5}{2}u_m(u_m-u^S)\times(-\frac{5}{2}u_m)(-\frac{3}{2}u_m-2u^S)\\[3mm]
&=\frac{25}{4}u_m^2[2u_m^2-4u_mu^S+2(u^S)^2+\frac{3}{2}u_m^2+\frac{1}{2}u_mu^S-2(u^S)^2]\\[3mm]
&=\frac{25}{4}u_m^2(\frac{7}{2}u_m^2-\frac{7}{2}u_m u^S)\\[3mm]
&=\frac{175}{8}u_m^4-\frac{175}{8}u_m^3u^S.
\end{aligned}
   \end{equation}  
When $\xi\in[\xi_1,\xi_*)$,  that is, $u^S\in[0, \frac{u_m}{2}),$ by \eqref{4.2},  
we have 
  \begin{equation}\label{a20}
 \begin{aligned}  
H_1&=\frac{5}{2}u_m\Big[\big(3u_m^2-3(u^S)^2\big)\big(-\frac{40}{u_m^2}(u^S)^3+\frac{30}{u_m}(u^S)^2-\frac{5}{2}u_m\big)\\[1mm]
 &\qquad\quad\  +3u^S\big(-\frac{10}{u_m^2}(u^S)^4+\frac{10}{u_m}(u^S)^3-\frac{5}{2}u_mu^{S}+\frac{5}{2}u_m^2\big)\\[1mm]
&\qquad\quad\  -\frac{1}{2}\big(-\frac{120}{u_m^2}(u^S)^2+\frac{60}{u_m}u^S\big)\big((u^S)^3-3u_m^2u^S+2u_m^3\big)\Big]\\[1mm]
&=\frac{375}{u_m}(u^S)^5-225(u^S)^4-750u_m(u^S)^3+750u_m^2(u^S)^2-\frac{525}{4}u_m^3u^S-\frac{75}{4}u_m^4,
  \end{aligned}
   \end{equation}  
   and 
     \begin{equation}\label{a21}
 \begin{aligned}
H_2&\di =w\Big[\big(-\frac{120}{u_m^2}(u^S)^2+\frac{60}{u_m}u^S\big)\big(u_m^2-\frac{3}{2}u_mu^S-(u^S)^2\big)-\frac{20}{u_m^2}(u^S)^4+\frac{20}{u_m}(u^S)^3\\
&\di \qquad\ \  -5u_mu^{S}+5u_m^2
+\big(-\frac{40}{u_m^2}(u^S)^3+\frac{30}{u_m}(u^S)^2-\frac{5}{2}u_m\big)\big(-\frac{3}{2}u_m-2u^S\big)\Big]\\
&\di =\Big[-\frac{10}{u_m^2}(u^S)^4+\frac{10}{u_m}(u^S)^3-\frac{5}{2}u_mu^{S}+\frac{5}{2}u_m^2\Big]\Big[\frac{180}{u_m^2}(u^S)^4+\frac{140}{u_m}(u^S)^3\\
&\di \qquad -255(u^S)^2+60u_mu^S+\frac{35}{4}u_m^2\Big]\\[1mm]
&\di =-\frac{1800}{u_m^4}(u^S)^8+\frac{400}{u_m^3}(u^S)^7+\frac{3950}{u_m^2}(u^S)^6-\frac{3600}{u_m}(u^S)^5\\[1mm]
&\di \quad+\frac{1225}{2}(u^S)^4+1075u_m(u^S)^3-\frac{1575}{2}u_m^2(u^S)^2+\frac{1025}{8}u_m^3u^S+\frac{175}{8}u_m^4.
  \end{aligned}
   \end{equation}  
   Thus, we get
   \begin{equation}\label{a22}
 \begin{aligned}
H_1+H_2=&\di -\frac{1800}{u_m^4}(u^S)^8+\frac{400}{u_m^3}(u^S)^7+\frac{3950}{u_m^2}(u^S)^6-\frac{3225}{u_m}(u^S)^5+\frac{775}{2}(u^S)^4\\[2mm]
&\di+325u_m(u^S)^3-\frac{75}{2}u_m^2(u^S)^2-\frac{25}{8}u_m^3u^S+\frac{25}{8}u_m^4.
  \end{aligned}
   \end{equation}  

{\bf Acknowledgment.}  The work of Feimin Huang was partially supported by National Key R$\&$D Program
of China No. 2021YFA1000800, and National Natural Sciences Foundation of China (NSFC) No. 12288201. The work of Yi Wang was partially supported by NSFC (Grant No. 12171459, 12288201, 12090014) and CAS Project for Young Scientists in Basic Research, Grant No. YSBR-031.

\
  
{\bf Conflict of Interest:} The authors declared that they have no conflicts of interest to this work.

\

{\bf Availability of data and material:} Data sharing not applicable to this article as no datasets were generated or analyzed during the current study.
     
{}

\begin{thebibliography}{xx}

\bibitem{BHR} M. Brio and J. K. Hunter. Rotationally invariant hyperbolic waves. Comm. Pure Appl. Math., 43 (8): 1037-1053, 1990.


\bibitem{KMY} K. Choi, M.-J. Kang, Y. Kwon, and A. Vasseur. Contraction for large perturbations of traveling waves in a hyperbolic-parabolic system arising from a chemotaxis model. Math. Models Methods Appl. Sci, 30(2):387-437, 2020.

\bibitem{LTP3} H. Freistühler and T.-P. Liu. Nonlinear stability of overcompresive shock waves in a rotationally invariant system of viscous conservation laws. Comm. Math. Phys., 153 (1): 147-158, 1993.

\bibitem{FS} H. Freist$\ddot{\mbox{u}}$hler, D. Serre,  $L^1$ stability of shock waves in scalar viscous conservation laws. Comm. Pure Appl. Math. 51 (1998), no. 3, 291-301.

\bibitem{HX} F. M. Huang and L. D.  Xu, Decay rate toward the traveling wave for scalar viscous conservation law, Commun. Math. Anal. Appl., 1 (2022), no. 3, 395-409.



\bibitem{IO} A. M. Il'in and O. A. Oleinik. Asymptotic behavior of solution  of the Cauchy problem for some quasilinear equations for large values of the time. Mat. Sb. (N.S.), 51 (93):191-216, 1960.


\bibitem{JGK} C. Jones, R. Gardner and T. Kapitula. Stability of traveling waves for non-convex scalar viscous conservation laws. Comm. Pure Appl. Math., 46:505-526, 1993.


\bibitem{VKM} M.-J. Kang and A. Vasseur. Contraction property for large perturbations of shocks of the barotropic Navier-Stokes system. J. Eur. Math. Soc., 23 (2):585-638, 2020. 

\bibitem{KV} M.-J. Kang and A. Vasseur.  $L^2$-contraction for shock waves of scalar viscous conservation laws. Ann. l'Institut Henri Poincar${\acute{\rm e}}$ C, Analyse non lineaire, 34 (1):139156, 2017.

\bibitem{MJ} M.-J. Kang.  $L^2$ -type contraction for shocks of scalar viscous conservation laws with strictly convex flux. J. Math. Pure Appl., 145:1-43, 2021.

\bibitem{KVW2} M.-J. Kang, A. Vasseur and Y. Wang. $L^2$-contraction of large planar shock waves for multi-dimensional scalar viscous conservation laws. Journal of Differential Equations, 267(5): 2737-2791, 2019.

\bibitem{KVW} M.-J. Kang, A. Vasseur and Y. Wang. Time-asymptotic stability of composite waves of viscous shock and rarefaction for barotropic Navier-Stokes equations. Adv. Math., 419:108963, 2023.

\bibitem{KVW1} M.-J. Kang, A. Vasseur, Y. Wang. Time-asymptotic stability of generic Riemann solutions for compressible Navier-Stokes-Fourier equations. arXiv:2306.05604 (2023).

\bibitem{MK} S. Kawashima and A. Matsumura. Stability of shock profiles in viscoelasticity with non-convex constitutive relations. Comm. Pure Appl. Math., 47 (12):1547-1569, 1994.


\bibitem{IPD} I. P. Lee-Bapty and D. G. Crighton. Nonlinear Wave Motion Governed by the Modified Burgers Equation. Philosophical Transactions of the Royal Society of London. Series A, Mathematical and Physical Sciences, 323:173-209, 1987.


\bibitem{LT} T.-P. Liu. Shock Waves. Graduate Studies in Mathematics.  American Mathematical Society, 215:106-109, 2021.


\bibitem{LHL} H. Liu. Asymptotic stability of shock profiles for non-convex convection-diffusion equation. Appl. Math. Letters, 10 (1):129-134, 1997.

\bibitem{M}
A.~Matsumura.
\newblock Waves in compressible fluids: viscous shock, rarefaction, and contact
  waves.
\newblock In {\em Handbook of mathematical analysis in mechanics of viscous
  fluids}, pages 2495--2548. Springer, Cham, 2018.

\bibitem{MN} A. Matsumura and K. Nishihara. Asymptotic stability of traveling waves for scalar viscous conservation laws with non-convex nonlinearity.  Comm. Math. Phys., 165:83-96, 1994.

\bibitem{MN1} A. Matsumura and K. Nishihara.  Asymptotics toward the rarefaction waves of the solutions of a one-dimensional model system for compressible viscous gas. Japan J. Appl. Math., 3:1-13, 1986.

\bibitem{MY} A. Matsumura and N. Yoshida.  Asymptotic behavior of solutions to the Cauchy problem for the scalar viscous conservation law with partially linearly degenerate flux. SIAM J. Math. Anal., 44 (4):2526-2544, 2012.

\bibitem{ME} M. Mei. Stability of shock profiles for non-convex scalar viscous conservation laws. Math. Mod. Meth. Appl. Sci., 5 (3):279-296,1995.



\bibitem{GAN} G. A. Nariboli and  W. C. Lin. A new type of Burgers' equation. ZAMM-J.  Appl. Math. Mech.
, 53 (8):505-510, 1973.

\bibitem{OR} S. Osher and J. Ralston. $ L_1$ stability of traveling waves with applications to convective porous media flow. Comm. Pure Appl. Math., 35:737- 751, 1982.



\bibitem{MEV} M. Teymur and E. Suhubi. Wave propagation in dissipative or dispersive non-linear media. IMA J. App. Math., 21(1):25-40, 1978.



\bibitem{We} H. F. Weinberger. Long-time behavior for a regularized scalar conservation law in the absence of genuine nonlinearity. Ann. Inst. Henri Poincar${\acute{\rm e}}$, 7:407-425, 1990.





\end{thebibliography}
\end{document}